\documentclass[12pt,reqno,oneside]{amsart}
\usepackage[letterpaper]{geometry}
\usepackage[rgb]{xcolor}
\usepackage[colorlinks, allcolors=blue]{hyperref}
\usepackage{amsmath,amsthm,amssymb,tikz,ifthen,enumitem,verbatim}

\allowdisplaybreaks % align/align* environment can break over pages
\setlist[itemize]{leftmargin=8mm}

\usepackage{cleveref} % Use \Cref to include "Theorem ", etc.
\newtheorem{thm}{Theorem}

\newtheorem{prop}[thm]{Proposition}
\newtheorem{conjecture}{Conjecture}
\theoremstyle{definition}

\newtheorem*{defn*}{Definition}
\newtheorem*{remark}{Remark}
\newcommand\<{\begin{equation}} \renewcommand\>{\end{equation}}

\usepackage[labelfont=sc, justification=centering, textformat=period]{caption}
\makeatletter \g@addto@macro\@floatboxreset\centering \makeatother % to center the figures themselves

\newcommand\abs[1]{\left\vert#1\right\vert}
\newcommand\C{{\mathbb C}}
\newcommand\D{{\mathbb D}}
\newcommand\R{{\mathbb R}}
\newcommand\Sb{{\mathbb S}}
\newcommand\Fc{{\mathcal F}}
\newcommand\Pc{{\mathcal P}}
\newcommand\Tc{{\mathcal T}}
\renewcommand\dots{...}
\newcommand\setbuilder[3][:]{\left\{\,#2 \,#1\, #3 \,\right\}}
\renewcommand\d{\mathrm{d}}

\newcommand\G{\Gamma}
\renewcommand\bar[1]{\,\overline{\!#1\!}\,}
\renewcommand\P{{\bar{P}}} \newcommand\Q{{\bar{Q}}} \newcommand\A{{\bar{A}}}
\renewcommand\i{\mathrm{i}}

\DeclareMathOperator\arccosh{arccosh}

\newcounter{commentcounter}\newcommand\COMMENT[2][red]{\stepcounter{commentcounter}\rlap{\smash{$^{\fcolorbox{#1}{#1!15}{\scriptsize\ifthenelse{\equal{#1}{green}}{\color{green!75!black}}{\color{#1}}\!\!{\bf\thecommentcounter}\!\!}}$}}\marginpar{\!\!\parbox{2.8cm}{\raggedright\small \ifthenelse{\equal{#1}{green}}{\color{green!67!black}}{\color{#1}}~\textbf{\thecommentcounter.}\,#2}}} \usepackage{silence} 

\begin{document}

\title[Rigidity and flexibility]{Rigidity and flexibility of entropies of boundary maps associated to Fuchsian groups}
\author{Adam Abrams}
\address{Faculty of Pure and Applied Mathematics, Wroc\l{}aw University of Science and Technology, Wroc\l{}aw 50370, Poland}
\email{the.adam.abrams@gmail.com}
\author{Svetlana Katok}
\address{Department of Mathematics, The Pennsylvania State University, University Park, PA 16802, USA}
\email{sxk37@psu.edu}
\author{Ilie Ugarcovici}
\address{Department of Mathematical Sciences, DePaul University, Chicago, IL 60614, USA}
\email{iugarcov@depaul.edu}

\thanks{The second author was partially supported by NSF grant DMS 1602409.}

\maketitle

\begin{abstract} Given a closed, orientable surface of constant negative curvature and genus $g \ge 2$, we study the topological entropy and measure-theoretic entropy (with respect to a smooth invariant measure) of generalized Bowen--Series boundary maps. Each such map is defined for a particular fundamental polygon for the surface and a particular multi-parameter.

We survey two strikingly different recent results by the authors: topological entropy is constant in this entire family (``rigidity''), while measure-theoretic entropy varies within Teichm\"uller space, taking all values (``flexibility'') between zero and a maximum, which is achieved on the surface that admits a regular fundamental $(8g-4)$-gon. We obtain explicit formulas for both entropies. The rigidity proof uses conjugation to maps of constant slope, while the flexibility proof---valid only for certain multi-parameters---uses the realization of geodesic flow as a special flow over the natural extension of the boundary map. We also present some new details pertaining to specific multi-parameters and specific polygons.
\end{abstract}

\section{Introduction}\label{sec intro}

Any closed, orientable, compact surface $S$ of constant negative curvature can be modeled as $S = \G\backslash\D$, where $\D = \{\, z \in \C : \abs z < 1 \,\}$ is the unit disk endowed with hyperbolic metric $2\abs{\mathrm{d}z}/(1-{\abs z}^2)$ and $\G$ is a finitely generated Fuchsian group of the first kind acting freely on $\D$ and isomorphic to $\pi_1(S)$.
Recall that geodesics in this model are half-circles or diameters orthogonal to $\Sb=\partial\D$, the circle at infinity. The geodesic flow $\widetilde\varphi^t$ on $\D$ is defined as an $\R$-action on the unit tangent bundle $T^1\D$ that moves a tangent vector along the geodesic defined by this vector with unit speed. The geodesic flow $\widetilde\varphi^t$ on $\D$ descends to the geodesic flow $\varphi^t$ on the factor $S=\G\backslash\D$ via the canonical projection of the unit tangent bundles. The orbits of the geodesic flow $\varphi^t$ are oriented geodesics on~$S$.

\begin{figure}[tb]
    \includegraphics[width=0.7\textwidth]{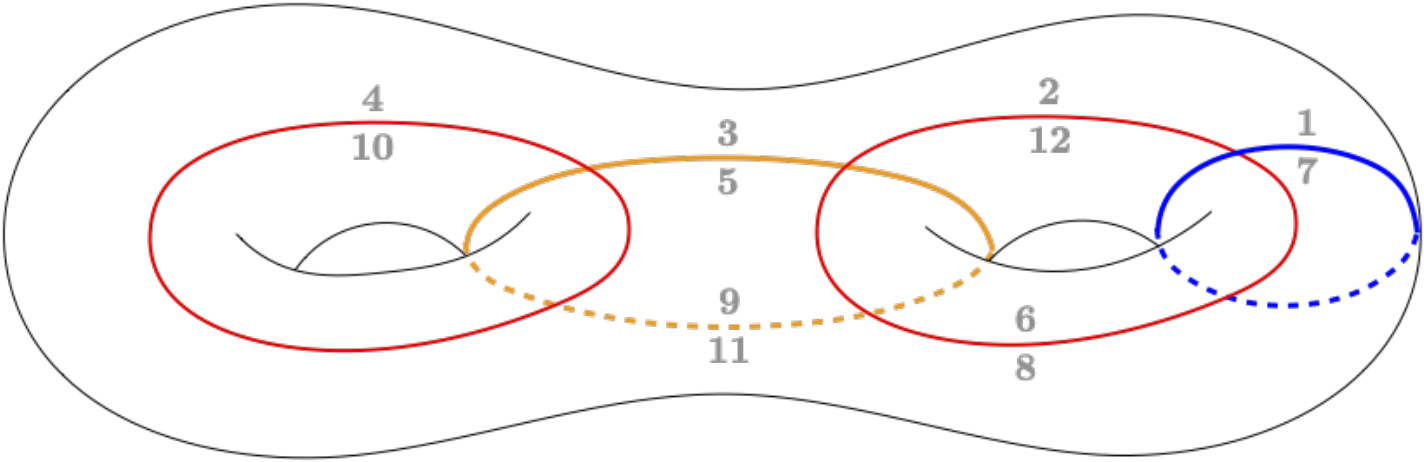} 
    \caption{Necklace of $2g$ geodesics on $S$ forming the sides of $\Fc$ for $g=2$}
    \label{fig necklace}
\end{figure}

A surface $S$ of genus $g \ge 2$ admits an  $(8g-4)$-sided fundamental polygon with a particular pairing of sides \eqref{sigma} obtained by cutting it with $2g$ closed geodesics that intersect in pairs ($g$ of them go around the ``holes'' of $S$ and another $g$ go around the ``waists''; see \Cref{fig necklace}). The existence of such a fundamental polygon $\Fc$ is an old result attributed to Dehn, Fenchel, Nielsen, and Koebe~\cite{W92,Ko29,BiS87}. Adler and Flatto~\cite[Appendix A]{AF91} give a careful proof of existence and special properties of $\Fc$ that we list below.

We label the sides of $\Fc$, which are geodesic segments, in a counterclockwise order by numbers $1 \le k \le 8g-4$ and label the vertices of $\Fc$ by $V_k$ so that side $k$ connects $V_k$ to $V_{k+1} \pmod {8g-4}$ (this gives us a {\em marking} of the polygon). 

We denote by $P_k$ and $Q_{k+1}$ the endpoints of the oriented infinite geodesic that extends side~$k$ to the circle at infinity~$\Sb$. The counter-clockwise order of endpoints on $\Sb$ is
\[ P_1, Q_1, P_2, Q_2, \dots, P_{8g-4}, Q_{8g-4}.
\label{eq: pq-partition} \]
The identification of the sides of $\Fc$ is given by the side pairing rule
\< \label{sigma} \sigma(k) := \left\{ \begin{array}{ll}
 4g-k \bmod (8g-4) & \text{ if $k$ is odd} \\
 2-k \bmod (8g-4) & \text{ if $k$ is even}.
\end{array} \right. \>
We denote by $T_k$ the M\"obius transformation pairing side~$k$ with side~$\sigma(k)$. The group~$\G$ is generated by $T_1,...,T_{8g-4}$.
\vspace*{-10pt}

\begin{figure}[htb]
	\includegraphics[height=3in]{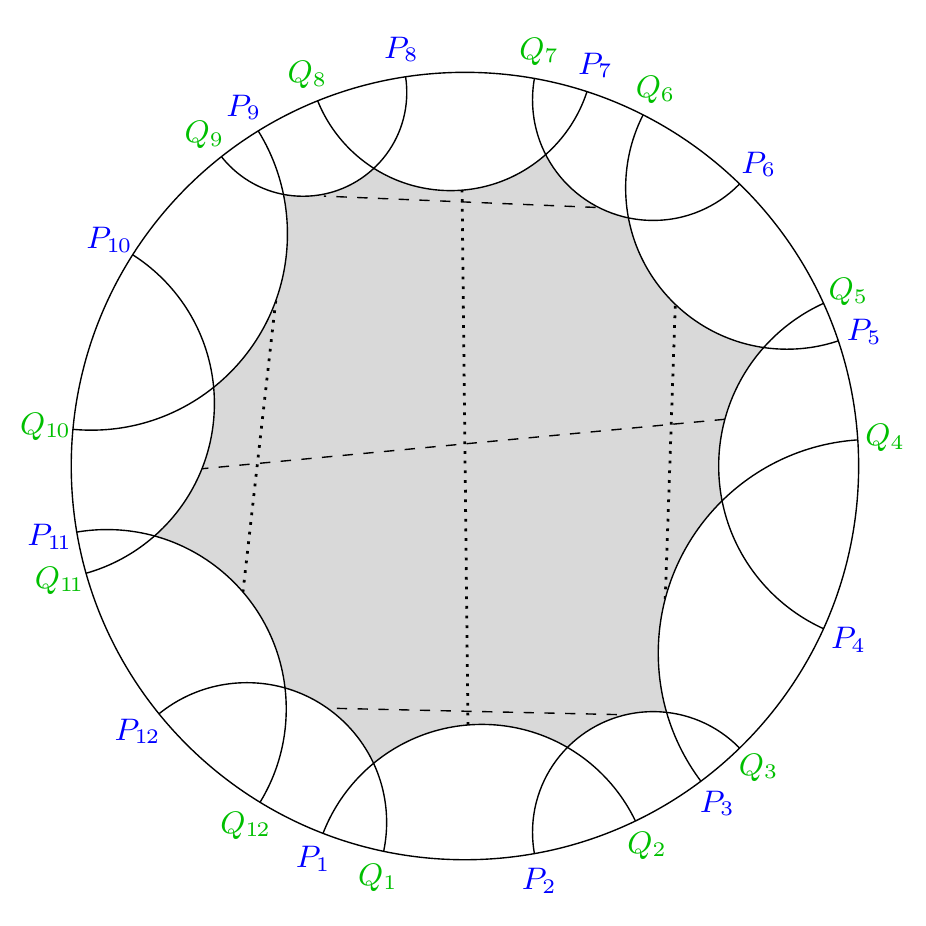} \quad
	\includegraphics[height=3in]{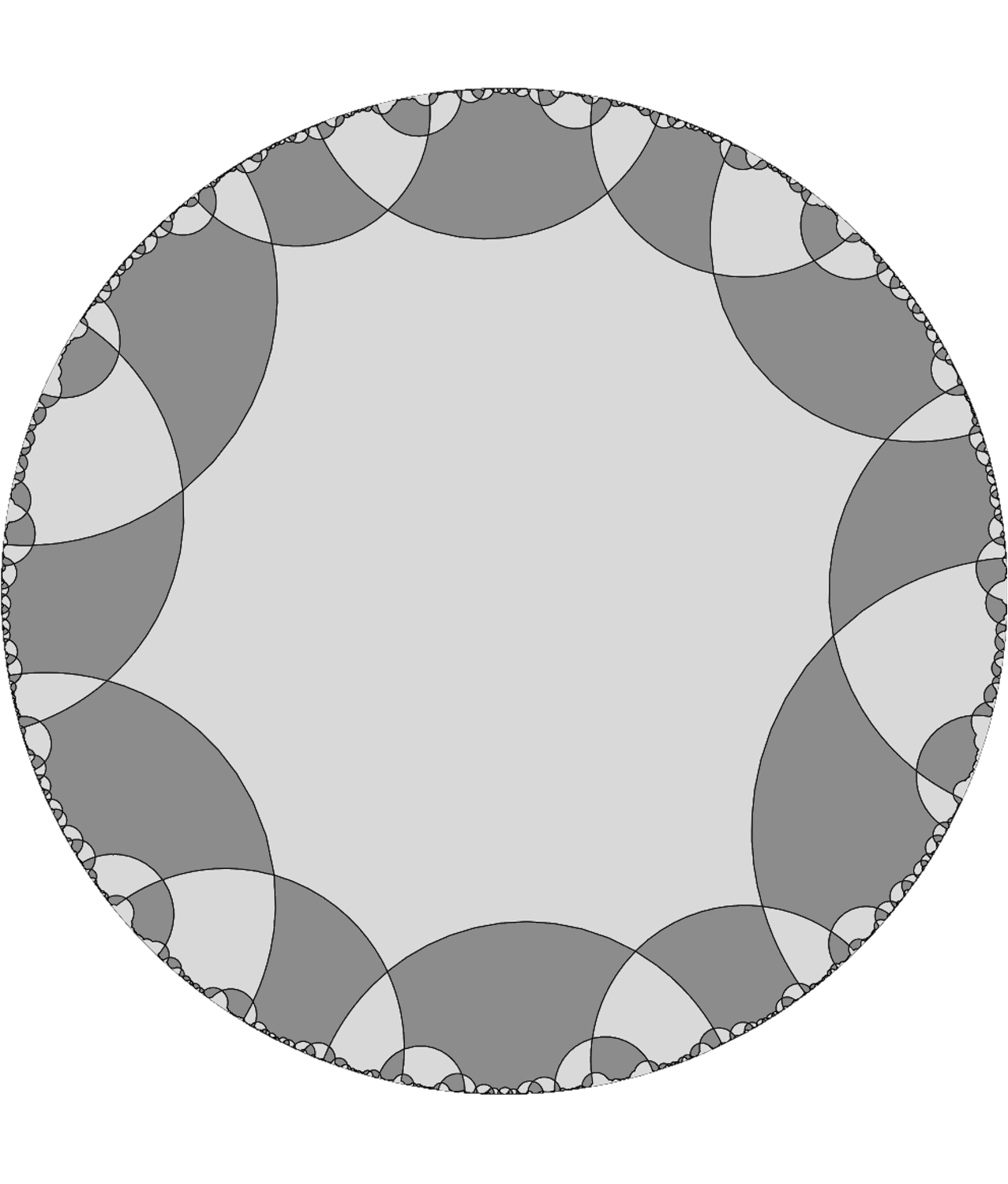}
	\vspace*{-10pt}
 \caption{An irregular polygon with dotted lines showing side identifications (left) and tessellation (right), genus $2$}
 \label{fig irregular sides}
\end{figure}

\pagebreak[4]
Note that $\Fc$ is not regular in general but satisfies the following properties:
\begin{enumerate}
\item the sides~$k$ and $\sigma(k)$  have equal length;
\item the angles at vertices~$k$ and~$\sigma(k)+1$ add up to $\pi$.
\end{enumerate}
We call fundamental $(8g-4)$-gons $\Fc$ satisfying properties (1) and (2) \emph{canonical}.
Property (2)  implies the ``extension condition,'' which is crucial for our analysis: the extensions of the sides of $\Fc$ do not intersect the interior of the tessellation $\gamma\Fc,\,\gamma\in \G$ (see \Cref{fig irregular sides}).

The polygon $\Fc$ is related to a regular $(8g-4)$-gon  $\Fc_{\rm reg}$ centered at the origin with all interior angles equal to $\frac\pi2$ (see~\cite[Figure~1]{AF91}) by the following construction. Let $S=\G\backslash \D$ be any compact surface of genus $g \ge 2$ and $S_{\rm reg}=\G_{\rm reg}\backslash \D$ be the surface of the same genus admitting the regular fundamental polygon $\Fc_{\rm reg}$. By the Fenchel--Nielsen Theorem~\cite{T72}, there exists an ori\-en\-ta\-tion-pre\-serv\-ing hom\-e\-o\-mor\-phism~$h$ from $\bar\D$ onto $\bar\D$ such that $\G=h\circ \G_{\rm reg} \circ h^{-1}$. The map $h\big|_\Sb$ is a homeomorphism of $\Sb$ preserving the order of the points $\P \cup \Q$, where $\P = \{ P_1,...,P_{8g-4}\}$ and $\Q = \{Q_1,...,Q_{8g-4}\}$.\footnote{In some works $h\big|_\Sb$ is also called a ``boundary map;'' this is different from our usage.} The sides of $\Fc$ are constructed by connecting $h(P_k)$ and $h(Q_{k+1})$ by geodesics, where $P_k$ and $Q_{k+1}$ correspond to~$\Fc_{\rm reg}$.
Conversely, for every marked canonical $(8g-4)$-gon and the associated Fuchsian group $\G$, by Euler's formula the genus of $\G\backslash \D$ is~$g$. 

This construction gives a less common representation of the Teichm\"uller space $\Tc(g)$ for $g \ge 2$: it is the space of marked canonical $(8g-4)$-gons in the unit disk~$\D$ up to an isometry of $\D$.\footnote{For a related representation of $\Tc(g)$ as the space of marked canonical $4g$-gons, see~\cite{S99}, following the earlier work~\cite{ZC80,T97}.} The topology on the space of polygons is that $\Pc_k \to \Pc$ if and only if the lengths of all sides converge and the measures of all angles converge. It is well-known~\cite{FM13,Funar} that $\dim(\Tc(g)) = 6g-6$; the following is a heuristic argument for this fact using $(8g-4)$-gons. The lengths of the identified pairs of sides are given by $4g-2$ real parameters, and $2g-1$ parameters represent the angles since four angles at each vertex are determined by one parameter. The dimension of the space $SU(1,1)$ of isometries of $\D$ is~$3$, so we have altogether $(4g-2)+(2g-1)-3=6g-6$ parameters.

\subsection*{The boundary map}

The main object of this paper is the multi-parameter family of generalized Bowen--Series boundary maps studied in~
\cite{KU17,AK19,A20,AKU-Flexibility,AKU-Rigidity}. For each $\Fc\in \Tc(g)$, 
recall that side~$k$ of~$\Fc$ is contained in the geodesic from $P_k$ to $Q_{k+1}$ and that $T_k$ maps side $k$ to side $\sigma(k)$; we define the \emph{boundary map} $f_\A:\Sb\to\Sb$
by 
\< \label{fA definition}
    f_\A (x)=T_k(x) \quad\text{if } x \in [A_k,A_{k+1}),
\>
where
\[ \A:=\{A_1,A_2,\dots,A_{8g-4}\}\quad\text{with}\quad A_k\in [P_k,Q_k]. \] 
When all $A_k = P_k$ we denote the map by $f_\P$ (see left of~\Cref{fig fP and CS} for its graph). The map $f_\P$ was extensively studied by Adler and Flatto in~\cite{AF91}; they call it ``the Bowen--Series boundary map,'' although Bowen and Series' construction~\cite{BS79} used $4g$-gons.

We study two dynamical invariants of~$f_\A$: the topological entropy and the measure-theoretic entropy with respect to a smooth (that is, Lebesgue-equivalent) measure. Our two main results are rigidity of the topological entropy (\Cref{thm rigidity}) and flexibility of the measure-theoretic entropy (\Cref{thm max and flexibility}). 

\begin{thm}[{Rigidity of topological entropy\,\cite{AKU-Rigidity}}] \label{thm rigidity} 
    Let $S=\G\backslash\D$ be a surface of genus $g \ge 2$ with fundamental polygon $\Fc \in \Tc(g)$.
    For any $\A = \{A_1,\dots,A_{8g-4}\}$ with $A_k \in [P_k,Q_k]$, the map $f_\A : \Sb \to \Sb$ has topological entropy \< \label{h top formula} h_\mathrm{top}(f_\A) = \log\!\big( 4g-3 + \sqrt{(4g-3)^2-1} \big). \>
\end{thm}

\begin{remark} The maps~$f_\A$ are not necessarily topologically conjugate since, according to~\cite{KU17}, the combinatorial structure of the orbits associated to the discontinuity points~$A_k$ can differ. \end{remark}

All $f_\A$ are piecewise continuous and piecewise monotone. Some maps in this family (such as those considered by Bowen and Series~\cite{BS79} and further studied by Adler and Flatto~\cite{AF91}) are Markov, and so the topological entropy can be calculated as the logarithm of the maximal eigenvalue of a transition matrix~\cite[Theorem 7.13]{W75} in these cases. Not all maps~$f_\A$ admit a Markov partition, but \eqref{h top formula} holds for all~$\A$ regardless. The proof of \Cref{thm rigidity} uses conjugation to maps of constant slope.

\bigskip
The flexibility result requires a restriction on the set of parameters, as defined below. Such restrictions are common in previous work on boundary maps:~\cite{AF91} uses only $\A=\P$ and $\A=\Q$,~\cite{A20} focuses on extremal parameters and their so-called duals, and~\cite{AK19} requires short cycles.

\begin{defn*} A multi-parameter $\A = \{A_1,\dots,A_{8g-4}\}$ is \emph{extremal} if for each $k$ either $A_k = P_k$ or $A_k = Q_k$. A multi-parameter $\A$ satisfies the \emph{cycle property} if each $A_k \in (P_k,Q_k)$ and there exist positive integers $m_k,n_k$ such that $f_\A^{m_k}(T_kA_k)=f_\A^{n_k}(T_{k-1} A_k)$. If $m_k=n_k=1$ then we say that $\A$ satisfies the \emph{short cycle property.}
\end{defn*}

Notice that the above properties are preserved by the Fenchel-Nielsen homeomorphism $h\big|_\Sb$ (see \cite{KU17e}), so it makes sense to talk about a multi-parameter $\A$ satisfying these properties in the Teichm\"uller space $\Tc(g)$.

If $\A$ is extremal or has short cycles, then $f_\A$ admits a unique smooth invariant ergodic measure $\mu_\A$ related to the Liouville measure for geodesic flow.
\begin{thm}[{Flexibility of measure-theoretic entropy\,\cite{AKU-Flexibility}}] \label{thm max and flexibility}~
    Let $S=\G\backslash\D$ be a surface of genus $g \ge 2$, and let~$\A$ be extremal or satisfy the short cycle property.\footnote{In~\cite{AKU-Flexibility} the main theorem is stated only for the classical Bowen--Series case $\A = \P$. As explained in \cite[Remark 5]{AKU-Flexibility}, the results also hold for $\A$ extremal or with short cycles, so we state the broader theorem here.}
    \begin{enumerate}[label=(\roman*)]
    \item For each $\Fc \in \Tc(g)$, the measure-theoretic entropy of $f_\A$ with respect to its smooth invariant probability measure $\mu_\A$ is 
    \< \label{h mu formula} h_{\mu_\A}(f_\A) = \frac{\pi^2(4g-4)}{\text{\small$\mathrm{Perimeter}$}(\Fc)} = \pi \cdot \frac{\mathrm{Area}(\Fc)}{\text{\small$\mathrm{Perimeter}$}(\Fc)}. \> 
	\item \label{item max} Among all $\Fc \in \Tc(g)$, the maximum value of $h_{\mu_\A}(f_\A)$ is achieved on the surface for which $\Fc$ is regular; this value is 
	\[ H(g) := h_{\mu_\A^\mathrm{reg}}(f_\A^\mathrm{reg})=\frac{\pi^2(4g - 4)}{(8g-4)\arccosh(1+2\cos\tfrac\pi{4g-2})}. \]
    \item \label{item flexibility} For any value $h \in (0,H(g)]$ there exists $\Fc \in \Tc(g)$ such that $h_{\mu_\A}(f_\A) = h$.
    \end{enumerate}
\end{thm}

\begin{remark} The formula~\eqref{h mu formula} shows that the measure-theoretic entropy remains constant under the change of the multi-parameter $\A$ within the considered class, so there is an aspect of rigidity to \Cref{thm max and flexibility} as well. We emphasize that $h_{\mu_\A}(f_\A)$ is flexible in the Teichm\"uller space. \end{remark}

The paper is organized as follows. Sections~\ref{sec rigidity} and \ref{sec flexibility} summarize the proofs of Theorems~\ref{thm rigidity} and~\ref{thm max and flexibility}, respectively. \Cref{sec two entropies} compares the two entropies (as was done in \cite{AKU-Flexibility}), and \Cref{sec sofic examples,sec genus 2} contain results about genus $2$ examples that are presented here for the first time. \Cref{sec open} lists some conjectures and open questions.

\section{Rigidity of topological entropy}\label{sec rigidity}

The notion of topological entropy was introduced by Adler, Konheim, and McAndrew in~\cite{AKM}. Their definition used covers and applied to compact Hausdorff spaces; Dinaburg~\cite{D70} and Bowen~\cite{B7173} gave definitions involving distance functions and separated sets, which are often more suitable for calculation. While these formulations of topological entropy were originally intended for continuous maps acting on compact spaces, Bowen's definition can actually be applied to piecewise continuous, piecewise monotone maps on an interval, as explained in~\cite{MZ}. The theory naturally extends to maps of the circle, where piecewise monotonicity is understood to mean local monotonicity or, equivalently, having a piecewise monotone lift to $\R$.

In~\cite{Parry64}, Parry introduced two probability measures for topological Markov chains: the first is the measure of maximal entropy and is commonly known as the ``Parry measure''; the second measure is not shift-invariant but is uniformly expanding on cylinders. In~\cite{Parry66}, Parry used this second measure to conjugate a piecewise monotone, strongly transitive interval map (not necessarily Markov) with positive topological entropy to a map with constant slope. Alsed\`a, Llibre, and Misiurewicz \cite{ALM,AM} generalized this construction to piecewise continuous, piecewise monotone interval maps, obtaining semi-conjugacy for non-transitive maps. We apply the results of~\cite{Parry66,AM} to the family of maps $f_\A$, defining the second of these measures, $\rho_\A$, in~\eqref{second Parry measure}.

Our original proof of \Cref{thm rigidity} in~\cite{AKU-Rigidity} uses several symmetry properties of the regular fundamental polygon $\Fc_{\rm reg}$ and the conjugacy $\psi_\P$ constructed for the boundary map $f_\P$ that is associated to the regular polygon. The main fact required for \Cref{thm rigidity} is \cite[Lemma~15(a)]{AKU-Rigidity}, and in hindsight this can be proved without appealing to any symmetry. Thus we can prove the rigidity result for any polygon $\Fc$ without using the Fenchel--Nielsen homeomorphism to relate the setup to the regular polygon. (Figures~\ref{fig fP and CS} and \ref{fig overlaps} were made using the regular polygon with genus $g=2$.)

\subsection*{Step 1: topological entropy for extremal parameters.} 
All $f_\A(P_k)$ and $f_\A(Q_k)$ belong to the set $\P \cup \Q$ (see~\cite[Proposition~2.2]{KU17}, originally~\cite[Theorem~3.4]{AF91}), so the partition of $\Sb$ into intervals $I_1,\dots,I_{16g-8}$ given by
\[ I_{2k-1} := [P_k, Q_k], \qquad I_{2k} := [Q_{k}, P_{k+1}], \qquad k=1,\dots,8g-4, \]
is a Markov partition for $f_\A$ for every extremal $\A$. Although they share a Markov partition, each extremal $\A$ has its own transition matrix $M_\A = (m_{i,j})$ given by
\< \label{transition matrix defn} m_{i,j} := \left\{\begin{array}{ll} 1 &\text{if } f_\A(I_i) \supset I_j \\ 0 &\text{otherwise.} \end{array}\right. \>
The transition matrices $M_\P$ and $M_\Q$ for genus 2 are shown in \Cref{fig PQ matrices}.
A word $\omega = (\omega_0,\dots,\omega_n)$ in the alphabet $\{1,\dots,16g-8\}$ is called \emph{$\A$-admissible} if $m_{\omega_i,\omega_{i+1}}=1$ for $i=0,\dots,n-1$.

\begin{figure}[htb]
\providecommand\0{{\color{black}0}}
\providecommand\1{{\color{black}\!\bf1\!}}

{\setlength\arraycolsep{1.7pt}
\renewcommand{\arraystretch}{0.9} %use for 12pt
\tiny$\left(\begin{array}{cccccccccccccccccccccccc}
\0 & \0 & \0 & \0 & \0 & \0 & \0 & \0 & \0 & \0 & \0 & \0 & \0 & \0 & \0 & \1 & \1 & \0 & \0 & \0 & \0 & \0 & \0 & \0 \\
\1 & \1 & \1 & \1 & \1 & \1 & \1 & \1 & \1 & \1 & \0 & \0 & \0 & \0 & \0 & \0 & \0 & \1 & \1 & \1 & \1 & \1 & \1 & \1 \\
\0 & \1 & \1 & \0 & \0 & \0 & \0 & \0 & \0 & \0 & \0 & \0 & \0 & \0 & \0 & \0 & \0 & \0 & \0 & \0 & \0 & \0 & \0 & \0 \\
\0 & \0 & \0 & \1 & \1 & \1 & \1 & \1 & \1 & \1 & \1 & \1 & \1 & \1 & \1 & \1 & \1 & \1 & \1 & \1 & \0 & \0 & \0 & \0 \\
\0 & \0 & \0 & \0 & \0 & \0 & \0 & \0 & \0 & \0 & \0 & \1 & \1 & \0 & \0 & \0 & \0 & \0 & \0 & \0 & \0 & \0 & \0 & \0 \\
\1 & \1 & \1 & \1 & \1 & \1 & \0 & \0 & \0 & \0 & \0 & \0 & \0 & \1 & \1 & \1 & \1 & \1 & \1 & \1 & \1 & \1 & \1 & \1 \\
\0 & \0 & \0 & \0 & \0 & \0 & \0 & \0 & \0 & \0 & \0 & \0 & \0 & \0 & \0 & \0 & \0 & \0 & \0 & \0 & \0 & \1 & \1 & \0 \\
\1 & \1 & \1 & \1 & \1 & \1 & \1 & \1 & \1 & \1 & \1 & \1 & \1 & \1 & \1 & \1 & \0 & \0 & \0 & \0 & \0 & \0 & \0 & \1 \\
\0 & \0 & \0 & \0 & \0 & \0 & \0 & \1 & \1 & \0 & \0 & \0 & \0 & \0 & \0 & \0 & \0 & \0 & \0 & \0 & \0 & \0 & \0 & \0 \\
\1 & \1 & \0 & \0 & \0 & \0 & \0 & \0 & \0 & \1 & \1 & \1 & \1 & \1 & \1 & \1 & \1 & \1 & \1 & \1 & \1 & \1 & \1 & \1 \\
\0 & \0 & \0 & \0 & \0 & \0 & \0 & \0 & \0 & \0 & \0 & \0 & \0 & \0 & \0 & \0 & \0 & \1 & \1 & \0 & \0 & \0 & \0 & \0 \\
\1 & \1 & \1 & \1 & \1 & \1 & \1 & \1 & \1 & \1 & \1 & \1 & \0 & \0 & \0 & \0 & \0 & \0 & \0 & \1 & \1 & \1 & \1 & \1 \\
\0 & \0 & \0 & \1 & \1 & \0 & \0 & \0 & \0 & \0 & \0 & \0 & \0 & \0 & \0 & \0 & \0 & \0 & \0 & \0 & \0 & \0 & \0 & \0 \\
\0 & \0 & \0 & \0 & \0 & \1 & \1 & \1 & \1 & \1 & \1 & \1 & \1 & \1 & \1 & \1 & \1 & \1 & \1 & \1 & \1 & \1 & \0 & \0 \\
\0 & \0 & \0 & \0 & \0 & \0 & \0 & \0 & \0 & \0 & \0 & \0 & \0 & \1 & \1 & \0 & \0 & \0 & \0 & \0 & \0 & \0 & \0 & \0 \\
\1 & \1 & \1 & \1 & \1 & \1 & \1 & \1 & \0 & \0 & \0 & \0 & \0 & \0 & \0 & \1 & \1 & \1 & \1 & \1 & \1 & \1 & \1 & \1 \\
\1 & \0 & \0 & \0 & \0 & \0 & \0 & \0 & \0 & \0 & \0 & \0 & \0 & \0 & \0 & \0 & \0 & \0 & \0 & \0 & \0 & \0 & \0 & \1 \\
\0 & \1 & \1 & \1 & \1 & \1 & \1 & \1 & \1 & \1 & \1 & \1 & \1 & \1 & \1 & \1 & \1 & \1 & \0 & \0 & \0 & \0 & \0 & \0 \\
\0 & \0 & \0 & \0 & \0 & \0 & \0 & \0 & \0 & \1 & \1 & \0 & \0 & \0 & \0 & \0 & \0 & \0 & \0 & \0 & \0 & \0 & \0 & \0 \\
\1 & \1 & \1 & \1 & \0 & \0 & \0 & \0 & \0 & \0 & \0 & \1 & \1 & \1 & \1 & \1 & \1 & \1 & \1 & \1 & \1 & \1 & \1 & \1 \\
\0 & \0 & \0 & \0 & \0 & \0 & \0 & \0 & \0 & \0 & \0 & \0 & \0 & \0 & \0 & \0 & \0 & \0 & \0 & \1 & \1 & \0 & \0 & \0 \\
\1 & \1 & \1 & \1 & \1 & \1 & \1 & \1 & \1 & \1 & \1 & \1 & \1 & \1 & \0 & \0 & \0 & \0 & \0 & \0 & \0 & \1 & \1 & \1 \\
\0 & \0 & \0 & \0 & \0 & \1 & \1 & \0 & \0 & \0 & \0 & \0 & \0 & \0 & \0 & \0 & \0 & \0 & \0 & \0 & \0 & \0 & \0 & \0 \\
\0 & \0 & \0 & \0 & \0 & \0 & \0 & \1 & \1 & \1 & \1 & \1 & \1 & \1 & \1 & \1 & \1 & \1 & \1 & \1 & \1 & \1 & \1 & \1
\end{array}\right)\quad
\left(\begin{array}{cccccccccccccccccccccccc}
\1 & \1 & \0 & \0 & \0 & \0 & \0 & \0 & \0 & \0 & \0 & \0 & \0 & \0 & \0 & \0 & \0 & \0 & \0 & \0 & \0 & \0 & \0 & \0 \\
\1 & \1 & \1 & \1 & \1 & \1 & \1 & \1 & \1 & \1 & \0 & \0 & \0 & \0 & \0 & \0 & \0 & \1 & \1 & \1 & \1 & \1 & \1 & \1 \\
\0 & \0 & \0 & \0 & \0 & \0 & \0 & \0 & \0 & \0 & \1 & \1 & \0 & \0 & \0 & \0 & \0 & \0 & \0 & \0 & \0 & \0 & \0 & \0 \\
\0 & \0 & \0 & \1 & \1 & \1 & \1 & \1 & \1 & \1 & \1 & \1 & \1 & \1 & \1 & \1 & \1 & \1 & \1 & \1 & \0 & \0 & \0 & \0 \\
\0 & \0 & \0 & \0 & \0 & \0 & \0 & \0 & \0 & \0 & \0 & \0 & \0 & \0 & \0 & \0 & \0 & \0 & \0 & \0 & \1 & \1 & \0 & \0 \\
\1 & \1 & \1 & \1 & \1 & \1 & \0 & \0 & \0 & \0 & \0 & \0 & \0 & \1 & \1 & \1 & \1 & \1 & \1 & \1 & \1 & \1 & \1 & \1 \\
\0 & \0 & \0 & \0 & \0 & \0 & \1 & \1 & \0 & \0 & \0 & \0 & \0 & \0 & \0 & \0 & \0 & \0 & \0 & \0 & \0 & \0 & \0 & \0 \\
\1 & \1 & \1 & \1 & \1 & \1 & \1 & \1 & \1 & \1 & \1 & \1 & \1 & \1 & \1 & \1 & \0 & \0 & \0 & \0 & \0 & \0 & \0 & \1 \\
\0 & \0 & \0 & \0 & \0 & \0 & \0 & \0 & \0 & \0 & \0 & \0 & \0 & \0 & \0 & \0 & \1 & \1 & \0 & \0 & \0 & \0 & \0 & \0 \\
\1 & \1 & \0 & \0 & \0 & \0 & \0 & \0 & \0 & \1 & \1 & \1 & \1 & \1 & \1 & \1 & \1 & \1 & \1 & \1 & \1 & \1 & \1 & \1 \\
\0 & \0 & \1 & \1 & \0 & \0 & \0 & \0 & \0 & \0 & \0 & \0 & \0 & \0 & \0 & \0 & \0 & \0 & \0 & \0 & \0 & \0 & \0 & \0 \\
\1 & \1 & \1 & \1 & \1 & \1 & \1 & \1 & \1 & \1 & \1 & \1 & \0 & \0 & \0 & \0 & \0 & \0 & \0 & \1 & \1 & \1 & \1 & \1 \\
\0 & \0 & \0 & \0 & \0 & \0 & \0 & \0 & \0 & \0 & \0 & \0 & \1 & \1 & \0 & \0 & \0 & \0 & \0 & \0 & \0 & \0 & \0 & \0 \\
\0 & \0 & \0 & \0 & \0 & \1 & \1 & \1 & \1 & \1 & \1 & \1 & \1 & \1 & \1 & \1 & \1 & \1 & \1 & \1 & \1 & \1 & \0 & \0 \\
\0 & \0 & \0 & \0 & \0 & \0 & \0 & \0 & \0 & \0 & \0 & \0 & \0 & \0 & \0 & \0 & \0 & \0 & \0 & \0 & \0 & \0 & \1 & \1 \\
\1 & \1 & \1 & \1 & \1 & \1 & \1 & \1 & \0 & \0 & \0 & \0 & \0 & \0 & \0 & \1 & \1 & \1 & \1 & \1 & \1 & \1 & \1 & \1 \\
\0 & \0 & \0 & \0 & \0 & \0 & \0 & \0 & \1 & \1 & \0 & \0 & \0 & \0 & \0 & \0 & \0 & \0 & \0 & \0 & \0 & \0 & \0 & \0 \\
\0 & \1 & \1 & \1 & \1 & \1 & \1 & \1 & \1 & \1 & \1 & \1 & \1 & \1 & \1 & \1 & \1 & \1 & \0 & \0 & \0 & \0 & \0 & \0 \\
\0 & \0 & \0 & \0 & \0 & \0 & \0 & \0 & \0 & \0 & \0 & \0 & \0 & \0 & \0 & \0 & \0 & \0 & \1 & \1 & \0 & \0 & \0 & \0 \\
\1 & \1 & \1 & \1 & \0 & \0 & \0 & \0 & \0 & \0 & \0 & \1 & \1 & \1 & \1 & \1 & \1 & \1 & \1 & \1 & \1 & \1 & \1 & \1 \\
\0 & \0 & \0 & \0 & \1 & \1 & \0 & \0 & \0 & \0 & \0 & \0 & \0 & \0 & \0 & \0 & \0 & \0 & \0 & \0 & \0 & \0 & \0 & \0 \\
\1 & \1 & \1 & \1 & \1 & \1 & \1 & \1 & \1 & \1 & \1 & \1 & \1 & \1 & \0 & \0 & \0 & \0 & \0 & \0 & \0 & \1 & \1 & \1 \\
\0 & \0 & \0 & \0 & \0 & \0 & \0 & \0 & \0 & \0 & \0 & \0 & \0 & \0 & \1 & \1 & \0 & \0 & \0 & \0 & \0 & \0 & \0 & \0 \\
\0 & \0 & \0 & \0 & \0 & \0 & \0 & \1 & \1 & \1 & \1 & \1 & \1 & \1 & \1 & \1 & \1 & \1 & \1 & \1 & \1 & \1 & \1 & \1
\end{array}\right)
$}
\caption{Transition matrices $M_\P$ (left) and $M_\Q$ (right) for $g=2$}
\label{fig PQ matrices}
\end{figure}

We compute the maximal eigenvalue of the transition matrices for all extremal parameters and hence the topological entropy for these Markov cases~\cite[Prop\-o\-si\-tion~7 and Cor\-ollary~8]{AKU-Rigidity}:
\begin{itemize}
    \item For any extremal multi-parameter $\A$, the number of $\A$-admissible words of length~$n$ grows as approximately $\lambda^n$, where \[ \lambda = 4g-3 + \sqrt{(4g-3)^2-1}, \] and therefore the maximal eigenvalue of~$M_\A$ is exactly $\lambda$.
    \item Thus $h_\text{top}(f_\A)=\log\lambda$ if $\A$ is extremal.
\end{itemize}

\subsection*{Step 2: conjugacy to a constant slope map.}

The following theorem combines several results of~\cite{Parry66,AM}, stated here for circle maps (as in~\cite{MSz}) instead of interval maps:
\begin{thm}\label{thm psi exists}
Given a piecewise monotone, piecewise continuous, topologically transitive map $f: \Sb \to \Sb$ of positive topological entropy $h>0$, there exists a unique (up to rotation of~$\Sb$) increasing homeomorphism $\psi:\Sb \to \Sb$ conjugating $f$ to a piecewise continuous map with constant slope $e^{h}$.
\end{thm}

\begin{figure}[hbt]
    \includegraphics[width=\textwidth]{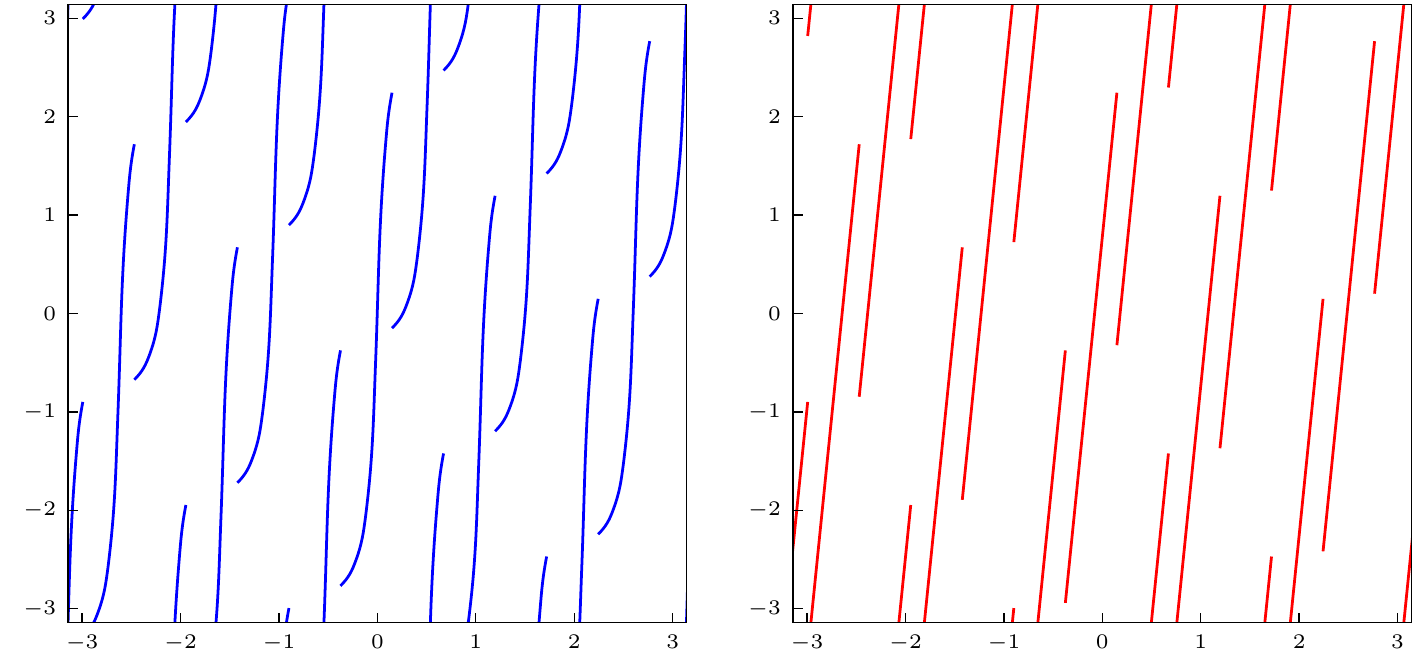}
    \caption{Plots of $f_\P(x)$ (left) and $\ell_\P(x)$ (right) for $g=2$}
    \label{fig fP and CS}
\end{figure}

Although for each parameter $\A$ (not necessarily extremal) the map $f_\A$ has, a priori, its own conjugacy and its own corresponding map of constant slope, we are especially interested in the cases $\A = \P$ and $\A = \Q$.
\begin{itemize}
\item The map $f_\P: \Sb \to \Sb$ is piecewise monotone, piecewise continuous, topologically transitive, and with positive topological entropy (see~\cite[Lemma 2.5]{BS79} and Step~1), so by \Cref{thm psi exists} there exists an increasing homeomorphism $\psi_\P : \Sb \to \Sb$ conjugating it to a map \[ \ell_\P := \psi_\P \circ f_\P \circ \psi_\P^{-1} \] with constant slope, see \Cref{fig fP and CS}. The map $\psi_\P$ is unique up to rotation of~$\Sb$, and the slope of $\ell_\P$ is exactly $\lambda = e^{h_\mathrm{top}(f_\P)}$.

\item The map $f_\Q: \Sb \to \Sb$ also satisfies the conditions of~\Cref{thm psi exists}, so there exists an increasing homeomorphism $\psi_\Q : \Sb \to \Sb$, unique up to rotation of~$\Sb$, conjugating it to a map~$\ell_\Q$ of constant slope. By Step~1, $\ell_\P$ and $\ell_\Q$ have the same slope.
\end{itemize}

When $f_\A$ is Markov---in particular when $\A$ is extremal---the construction of the conjugacy $\psi_\A$ follows the classical work of Parry~\cite{Parry66}, also used in the proof of~\cite[Lem\-ma~5.1]{AM}: 
let $\lambda,v$ be the maximal eigenpair for the transition matrix $M_\A$ and define the measure~$\rho_\A$ on the shift space~$X_\A$ for any non-empty cylinder~$C_\A^\omega$ (see~\cite[Section 4]{AKU-Rigidity}) as
\< \label{second Parry measure}
	\rho_\A\big(C_\A^{(\omega_0,\dots,\omega_n)}\big) = \frac{v_{\omega_n}}{\lambda^n}. 
\>
One defines the push-forward measure~$\rho'_\A$ on $\Sb=[-\pi,\pi]$ as
\[ \label{mu'} \rho'_\A(E) = \rho_\A\big( \phi^{-1}_\A(E) \big) \qquad\text{for Borel $E$}, \]
where $\phi_\A: X_\A \to \Sb$ is the symbolic coding map, that is, $\phi_\A(\omega) = \bigcap_{i=0}^\infty f_\A^{-i}(I_{\omega_i})$.
The conjugacy $\psi_\A:\Sb \to \Sb$ in Markov cases is then given by 
\< \label{psi formula from Parry} 
\psi_\A(x) := -\pi + 2 \pi \cdot \rho'_\A\big([-\pi,x]\big).
\>

\smallskip
It turns out that the maps $\psi_\P$ and $\psi_\Q$ thus constructed coincide:
\begin{thm}[{\cite[Theorem~12]{AKU-Rigidity}}] \label{thm psi P Q are equal} For all $x\in\Sb$, $\psi_\P(x) = \psi_\Q(x)$. \end{thm}
This is the most technically difficult part of the argument. The proof is presented in~\cite[Appendix~A]{AKU-Rigidity}. We also prove many symmetric properties of $\psi_\P$ (\cite[Propositions~14 and~15]{AKU-Rigidity}) that are visually evident in \cite[Figure~4]{AKU-Rigidity}.

\subsection*{Step 3: completion of the proof.}

\begin{figure}[hbt]
 \includegraphics[width=0.7\textwidth]{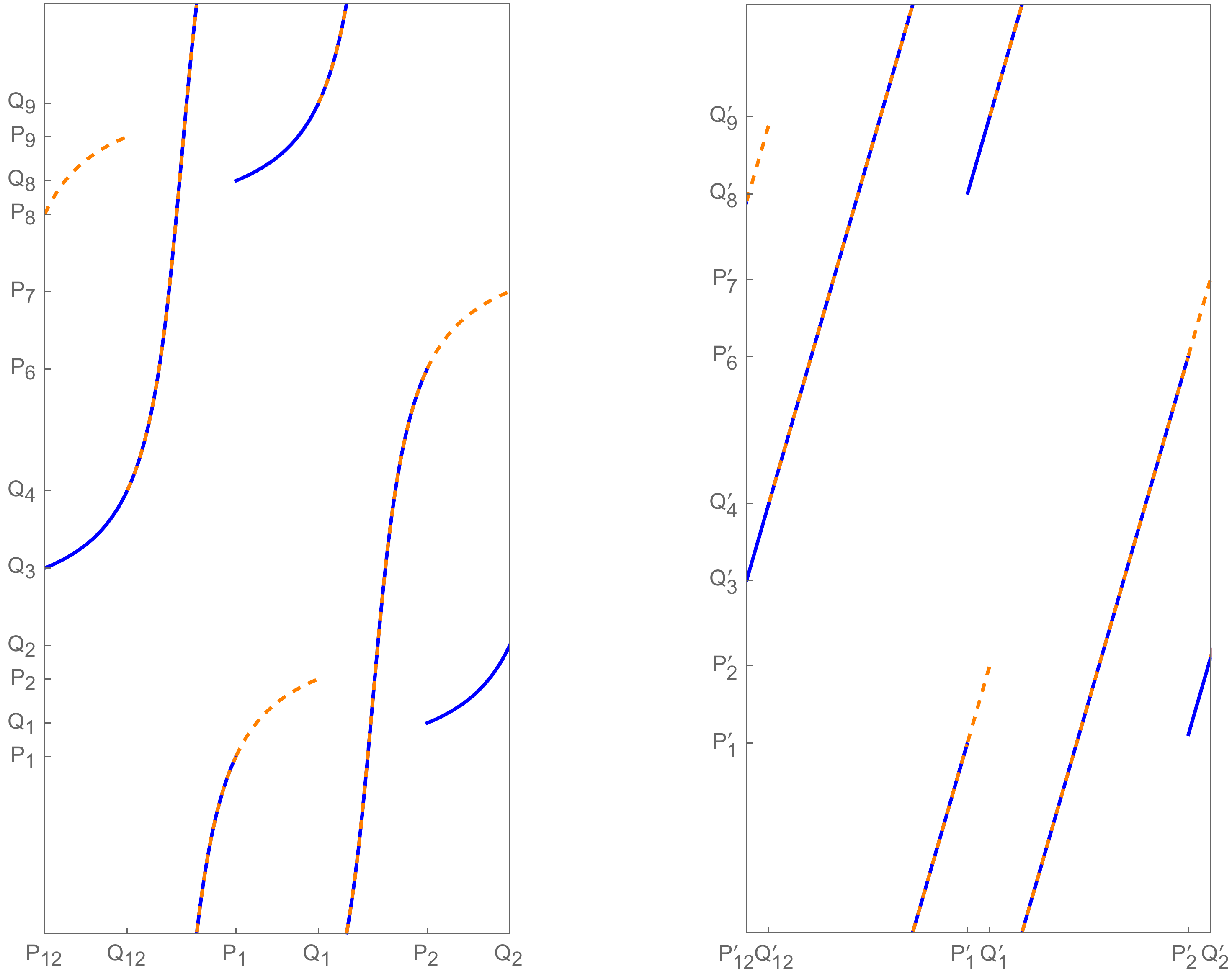}
 \caption{Left: part of the graphs of $f_\P$ (blue) and $f_\Q$ (dashed orange). Right: $\ell_\P$ (blue) and $\ell_\Q$ (dashed orange)}
 \label{fig overlaps}
\end{figure}

Denote $P'_k = \psi_\P(P_k)$ and $Q'_k = \psi_\P(Q_k)$.
By construction, the map $\psi_\P \circ T_k \circ \psi_\P^{-1}$ is linear (with slope $\lambda$) on $[P'_k,P'_{k+1})$ since $[P_k,P_{k+1})$ is the interval where $f_\P$ acts by $T_k$. The crucial fact that is necessary for our proof of rigidity is that $\psi_\P \circ T_k \circ \psi_\P^{-1}$ is linear on the longer interval
\[ [P'_k,Q'_{k+1}) \supset [P'_k,P'_{k+1}).  \]
To see this, note that $\psi_\Q \circ T_k \circ \psi_\Q^{-1}$ is linear on $[Q'_k,Q'_{k+1})$ by construction and that, by \Cref{thm psi P Q are equal}, $\psi_\P =  \psi_\Q$. Thus $\psi_\P \circ T_k \circ \psi_\P^{-1} = \psi_\Q \circ T_k \circ \psi_\Q^{-1}$ is linear on $[P'_k,P'_{k+1}] \cup [Q'_k,Q'_{k+1}] = [P'_k,Q'_{k+1}]$. (See \Cref{fig overlaps}.)

\begin{proof}[Proof of \Cref{thm rigidity}]
Let $\A = \{A_1,...,A_{8g-4}\}$ with each $A_k \in [P_k,Q_k]$ be arbitrary. For each $k$, the map \[ \psi_\P \circ f_\A \circ \psi_P^{-1} \] (note the use of $\psi_\P$ with $f_\A$) is linear with slope $\lambda$ on $[A'_k,A'_{k+1})$ because $\psi_\P \circ T_k \circ \psi_\P^{-1}$ is linear on $[A'_k,A'_{k+1}) \subset [P'_k,Q'_{k+1})$. Thus on whole full set $\bigcup_{k=1}^{8g-6} [A'_k,A'_{k+1}) = \Sb$ the map $\psi_\P \circ f_\A \circ \psi_\P^{-1}$ is linear with slope $\lambda$, and so, by~\cite[Theorem~{$3'$}]{MSz}, the topological entropy of $f_\A$ is $\log\lambda$. \end{proof}

\section{Flexibility of measure-theoretic entropy}\label{sec flexibility}
A few years ago, Anatole Katok suggested a new area of research---or, at the very least, a new viewpoint---called the ``flexibility program,'' which can be broadly formulated as follows: under properly understood general restrictions, within a fixed class of smooth dynamical systems, some dynamical invariants take arbitrary values. Taking this point of view, it is natural to ask how the mea\-sure-theo\-re\-tic entropy $h_{\mu_\P}(f_\P)$ changes in $\Tc(g)$. \Cref{thm max and flexibility} addresses this question, and we sketch its proof in this section.

\subsection*{Step 1: the smooth invariant measure.} If $\A$ is extremal (e.g., $\A=\P$) or has short cycles, then the smooth invariant measure $\mu_\A$ for $f_\A$ can be described as two-step projection of the Liouville measure for the geodesic flow~\cite{AK19,A20}.

\begin{itemize}
	\item On the unit tangent bundle $T^1\D$, parameterized by $(x,y,\theta)$ with $x+\i y \in \D$ the base-point of tangent vector and $\theta$ its angle with the real axis, the Liouville volume \[ \d \omega = \frac{4\,\d x\,\d y\,\d\theta}{(1-x^2-y^2)^2} \] comes from the hyperbolic measure on $\D$ and is invariant under geodesic flow.
	
	\item The measure \[ \d\nu = \dfrac{\abs{\d u}\abs{\d w}}{\abs{u-w}^2} \] is a smooth measure on the space of oriented geodesics on~$\D$ (modeled as $\{(u,w)\in\Sb\times\Sb : u \ne w\}$) and is preserved by M\"obius transformations.\footnote{The formula $\d\nu = |\d u||\d w|/|u-w|^2$ is used when $u,w \in \{z \in \C:|z|=1\}$. If we consider $u,w \in (-\pi,\pi]$ as arguments, then $\d\nu = \d u\d w/(2-2\cos(u-w))$.} Sometimes called ``geodesic current,'' this measure was most probably first considered by E.~Hopf~\cite{H36} and was later used by Sullivan~\cite{S79}, Bonahon~\cite{B88}, Adler--Flatto~\cite{AF91}, and the current authors~\cite{KU17,AK19}.
	
	\item Using coordinates $(u,w,s)$ on $T^1\D$, where $s$ is arclength along a geodesic, instead of coordinates $(x,y,\theta)$, the measure
	\[ \d m = \d\nu\,\d s \]
	is also preserved by geodesic flow and is a multiple of the Liouville volume~$\d \omega$. Specifically, $\d m = \tfrac12 \d \omega$ (see~\cite[Appendix~A2]{B88}).\footnote{\,The constant relating $d\omega$ and $\d m$ was given incorrectly as $1/4$ in~\cite[page 250]{AF91}. Following that, $1/4$ was used in~\cite[Proposition~10.1]{AK19}.} %In the setting of variable negative curvature surfaces, $\d m$ is the ``Patterson--Sullivan'' measure and $\d\omega$ is the ``Bowen--Margulis'' measure; for constant curvature these measures coincide up to a scalar multiple.
\end{itemize}

\smallskip
Generalizing Adler and Flatto's ``rectilinear map'' from~\cite{AF91}, we define the map $F_\A:\Sb\times\Sb\setminus\Delta \to \Sb\times\Sb\setminus\Delta$, where $\Delta$ is the diagonal $\setbuilder{(w,w)}{w\in\Sb}$, by
\< \label{FP defn} F_\A(u,w) = (T_ku,T_kw) \qquad\text{if }w \in [A_k,A_{k+1}). \>
In~\cite{KU17}, the authors showed that $F_\A$ admits a global attractor $\Omega_\A$ with finite rectangular structure if $\A = \P$ or if $\A$ satisfies the short cycle property. Adler and Flatto had previously shown that~$F_\P$ has as an invariant domain with finite rectangular structure (which we call~$\Omega_\P$, see \Cref{fig two domains}), and this was extended to $F_\A$ for all extremal parameters in~\cite{A20}. The re\-stric\-tion of $F_\A$ to ${\Omega_\A}$ is the natural extension map of~$f_\A$; we will often denote $F_\A\big|_{\Omega_\A}$ by simply $F_\A$.
\begin{itemize}
	\item The normalized measure \[ \label{dnuP} \d\nu_\A := \frac{\d\nu}{\int_{\Omega_\A}\d\nu} \] is by construction a smooth invariant probability measure for $F_\A$.
\end{itemize}
Because $f_\A$ is a factor of $F_\A$ (projecting on the second coordinate), its smooth invariant probability measure $\mu_\A$ is obtained as a projection of $\nu_\A$.

\subsection*{Step 2: formula for the entropy.} 
The geodesic flow on $S$ can be realized as a special flow over a cross-section that is parametrized by $\Omega_\A$, and the first return map to this cross-section acts exactly as $F_\A : \Omega_\A \to \Omega_\A$ (see \cite[Section~4]{AK19} for short cycles and \cite[Section~2.2]{A20} for extremal).
Using this realization along with Abramov's formula and the Ambrose--Kakutani theorem, we have from~\cite[Proposition~10.1]{AK19} (with a corrected constant) that 
\[ 
h_{\nu_\A}(F_\A) = \frac{\pi^2(4g-4)}{\int_{\Omega_\A}\d\nu}.
\]
Since $F_\A$ is the natural extension of $f_\A$, the entropies $h_{\nu_\A}(F_\A)$ and $h_{\mu_\A}(f_\A)$ are equal, and since $\mathrm{Area}(\Fc)=2\pi(2g-2)$ by the Gauss--Bonnet formula, we have
\< \label{h from AK} h_{\mu_\A}(f_\A) = \pi \cdot \frac{\mathrm{Area}(\Fc)}{\int_{\Omega_\A}\d\nu}. \>
\begin{figure}[tb]
	\includegraphics[width=0.67\textwidth]{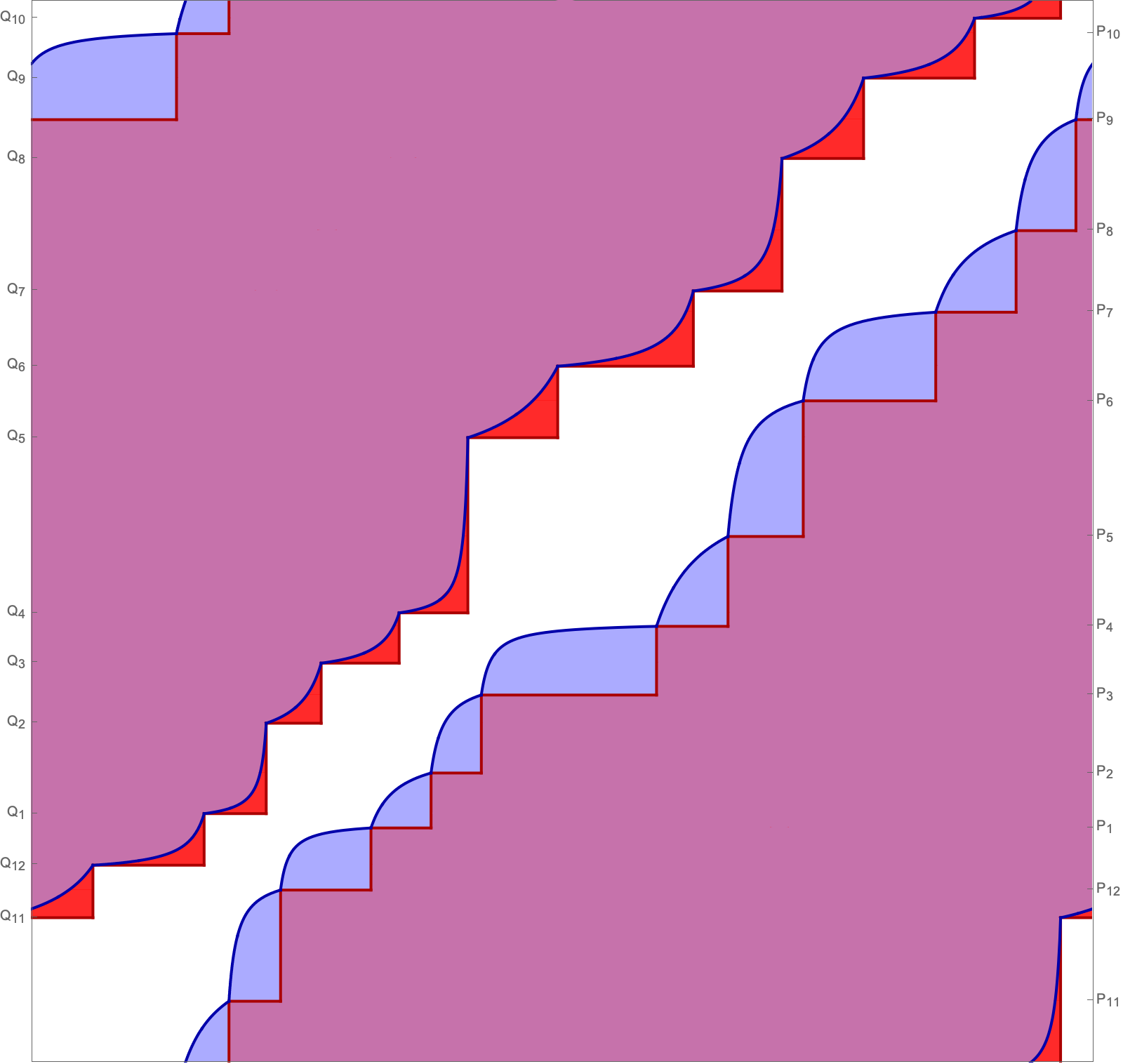}
 \caption{Arithmetic set $\Omega_\P$ in red-and-purple, geometric set $\Omega_\mathrm{geo}$ in blue-and-purple (for an irregular polygon with $g=2$)}
 \label{fig two domains}
\end{figure}

\noindent The integral $\int_{\Omega_\A}\d\nu$ can be explicitly computed using the finite rectangular structure of $\Omega_\A$, but this does not easily relate to the (hyperbolic) perimeter of $\Fc$.

Adler--Flatto~\cite{AF91} introduced another map, called the ``curvilinear map'' (or ``geometric map'' in~\cite{AK19}): denoting by $uw$ the geodesic from $u$ to $w$, the map is defined on the set
\[ \Omega_{\rm geo}:=\setbuilder{ (u,w) }{ uw \text{ intersects }\Fc} \;\;\subset\;\; \Sb\times\Sb\setminus\Delta \]
(shown in \Cref{fig two domains}) and is given by
\[ F_{\rm geo}(u,w) = (T_ku,T_kw) \quad\text{if $uw$ exits $\Fc$ through side~$k$.} \]

There is a bijection $\Phi_\A: \Omega_{\rm geo} \to \Omega_\A$ that acts piecewise by M\"obius transformations (see~\cite[Proposition 3.4]{AK19} for short cycles and~\cite[Proposition~12]{A20} for extremal) and since M\"obius transformations preserve the measure $\nu$, it follows that 
\< \label{PG} 
    \int_{\Omega_\A}\d\nu = \int_{\Omega_{\rm geo}}\d\nu.
\>

Given \eqref{PG}, we now want to show that $\int_{\Omega_{\rm geo}}\d\nu$ is equal to the (hyperbolic) perimeter of~$\Fc$. We use the following fact proved by F.~Bonahon~\cite[Appendix~A3]{B88}:
\begin{itemize} 
	\item For any oriented geodesic segment $s$ on $\D$, 
	\[ \int_{\Psi^+\!(s)}\d\nu = \mathrm{length}(s), \]
	where $\Psi^+\!(s)$ is the set of oriented geodesics intersecting $s$ with the oriented angle at the intersection between $0$ and $\pi$.
\end{itemize}
The domain $\Omega_{\rm geo}$ of the geometric map $F_{\rm geo}$ can be decomposed as $\Omega_{\rm geo} = \bigcup_{k=1}^{8g-4} {\mathcal G}_k$, where \[ {\mathcal G}_k = \setbuilder{ (u,w) }{ uw\text{ exits $\Fc$ through side }k } = \Psi_+(\text{side }k) \]
(these ``strips'' are shown in~\cite[Figure~3]{AK19}). Thus from Bohanon's result we immediately get 
\[ \label{Gperim} \int_{\Omega_{\rm geo}}\d\nu = \sum_{k=1}^{8g-4} \int_{{\mathcal G}_k}\d\nu = \sum_{k=1}^{8g-4} \mathrm{length}(\text{side }k) = \mathrm{Perimeter}(\Fc). \]
Combining this with \eqref{PG}, one can replace $\int_{\Omega_\A}\d\nu$ by the perimeter of $\Fc$ in the denominator of \eqref{h from AK}, and we obtain the formula \eqref{h mu formula}, that is,
\[
 	h_{\mu_\A}(f_\A) 
 	= \frac{\pi^2(4g-4)}{\text{\small$\mathrm{Perimeter}$}(\Fc)}
 	= \pi \cdot \frac{\mathrm{Area}(\Fc)}{\text{\small$\mathrm{Perimeter}$}(\Fc)}.
 \]

\subsection*{Step 3: maximum of the entropy.}
To prove that \Cref{thm max and flexibility}\ref{item max} follows from
\eqref{h mu formula} we only need to show that for each genus $g$ the perimeter of $\Fc$ in $\Tc(g)$ is minimized on the regular polygon.
\begin{itemize}
    \item Isoareal Inequality: {\em among all hyperbolic polygons with a given area and number of sides, the regular polygon has the smallest perimeter.}
    More precisely, Ku-Ku-Zhang prove~\cite[Theorem 1.2(a)]{Ku} that for a hyperbolic $n$-gon $\Pc_n$,
    \[
    \mathrm{Perimeter}(\Pc_n)^2\ge 4d_n\mathrm{Area}(\Pc_n), \text{ where }d_n=n\tan\left(\frac{\mathrm{Area}(\Pc_n)}{2n}\right),
    \]
    with equality achieved on a regular polygon.
    The Isoareal Inequality follows immediately: $\mathrm{Area}(\Pc_n)$ and $n$ are constant, so the right-hand side $4d_n\mathrm{Area}(\Pc_n)$ is constant and thus the perimeter is minimized when $\Pc_n$ is regular.
    \item In our setting, $\Fc=\Pc_n$ with $n=8g-4$, and $\mathrm{Area}(\Fc)=2\pi(2g-2)$ is constant in $\Tc(g)$, so the Isoareal Inequality implies that the perimeter of $\Fc$ is minimized when $\Fc$ is regular.
\end{itemize}  
The expression for the maximum value $H(g)$ in \Cref{thm max and flexibility} comes directly from \eqref{h mu formula}, with
\[ \arccosh\!\big(1+2\cos\tfrac{\pi}{4g-2}\big) \]
being the length of a single side of the regular $(8g-4)$-gon. This completes the proof of \Cref{thm max and flexibility}\ref{item max}.

\subsection*{Step 4: flexibility of the entropy.}

\begin{figure}[tb]
    \includegraphics[width=0.67\textwidth]{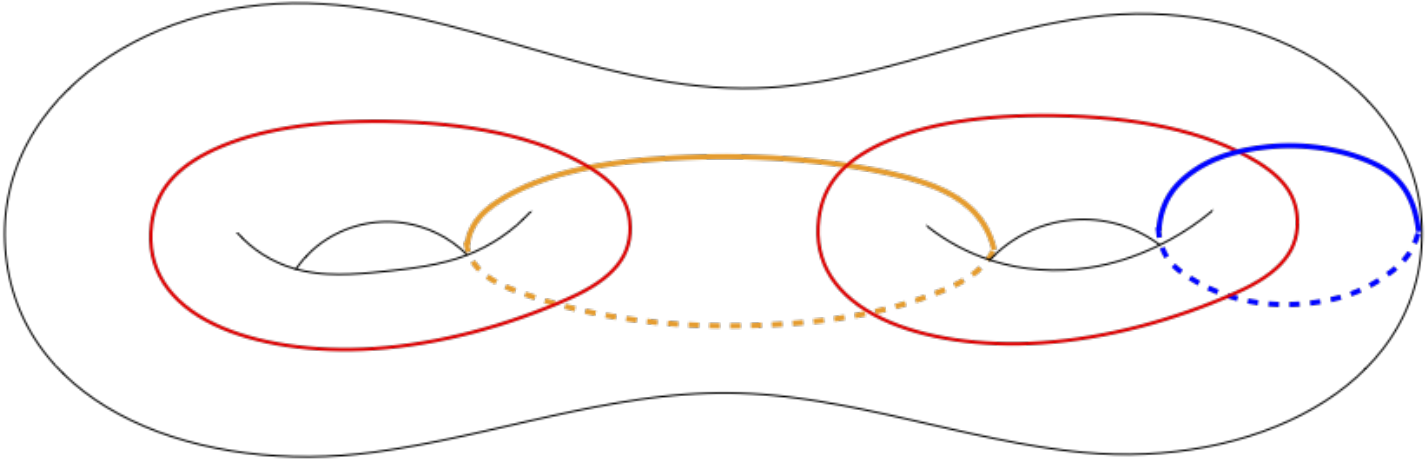} \\[0.25em]
    \includegraphics[width=0.67\textwidth]{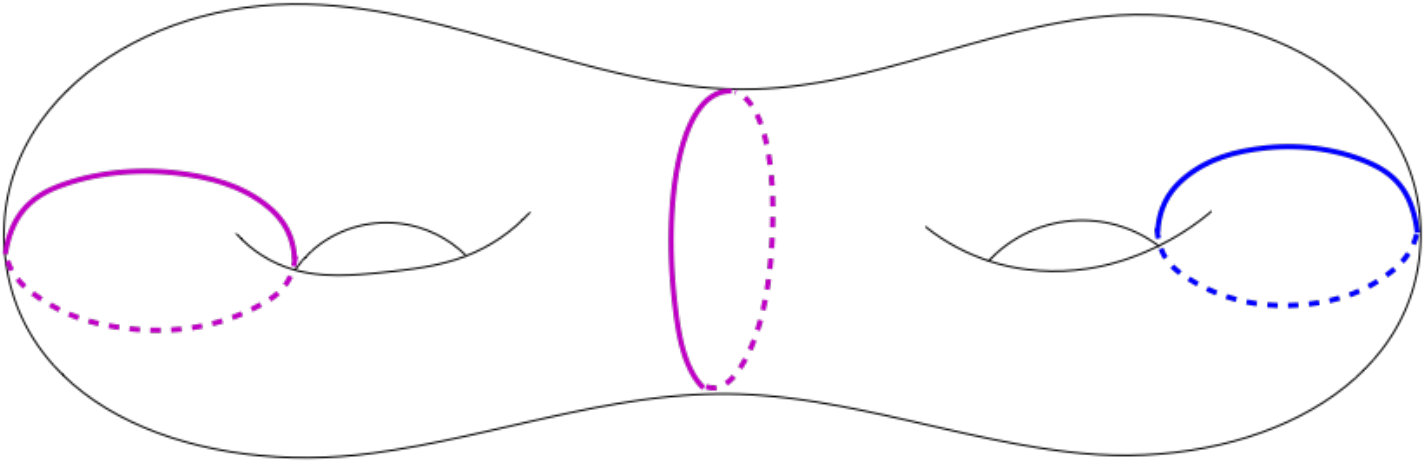}
    \caption{Necklace of $2g$ geodesics on $S$ forming the sides of $\Fc$ (top) and decomposition of $S$ into $2g-2$ pairs of pants by $3g-3$ non-intersecting geodesics (bottom) for $g=2$}
    \label{fig loops}
\end{figure}

The space~$\Tc(g)$ is ho\-me\-o\-morph\-ic to $\R^{6g-6}$, and a standard way to parametrize $\Tc(g)$ is through Fenchel--Nielsen coordinates (see the classical manuscript recently published in~\cite{FN}). The surface $S$ can be decomposed into $2g-2$ pairs of pants by $3g-3$ non-intersecting closed geodesics; this decomposition is shown for $g=2$ in the bottom of \Cref{fig loops} (for genus $3$, see~\cite[Figure~4]{AKU-Flexibility}). The lengths of these geodesics can be manipulated independently (they form $3g-3$ of the $6g-6$ coordinates) and can take arbitrarily large values. We take one of these geodesics to also be a geodesic from the necklace described in \Cref{sec intro} that corresponds to one entire side of $\Fc$ (this shared geodesic is on the far right in both parts of \Cref{fig loops}). Since the length of this side---one of the Fenchel--Nielsen coordinates---can be made arbitrarily large, the perimeter of $\Fc$ can also be made arbitrarily large, which by~\eqref{h mu formula} means that $h_{\mu_\A}(f_\A)$ can be made arbitrarily small.

Using the topology on the Teichm\"uller space~$\Tc(g)$ as the space of marked canonical $(8g-4$)-gons, we see that the perimeter of~$\Fc$ varies continuously within~$\Tc(g)$. From~\eqref{h mu formula} we conclude the continuity of the entropy $h_{\mu_\A}(f_\A)$ within~$\Tc(g)$. By the Intermediate Value Theorem, $h_{\mu_\A}(f_\A)$ must take on all values between $0$ and its maximum; this is precisely the claim of \Cref{thm max and flexibility}\ref{item flexibility}.

\section{Additional results} \label{sec new}

\subsection{Comparison of entropies}\label{sec two entropies}

By the Variational Principle, the measure-the\-or\-e\-tic entropy of a map can never exceed its topological entropy. For some classical systems, the Lebesgue measure is both a smooth invariant measure and also the measure of maximal entropy, but the boundary maps $f_\A: \Sb \to \Sb$ provide examples where the ergodic smooth invariant probability measure is not the measure of maximal entropy. 

\begin{figure}[bth]
 \includegraphics[width=0.8\textwidth]{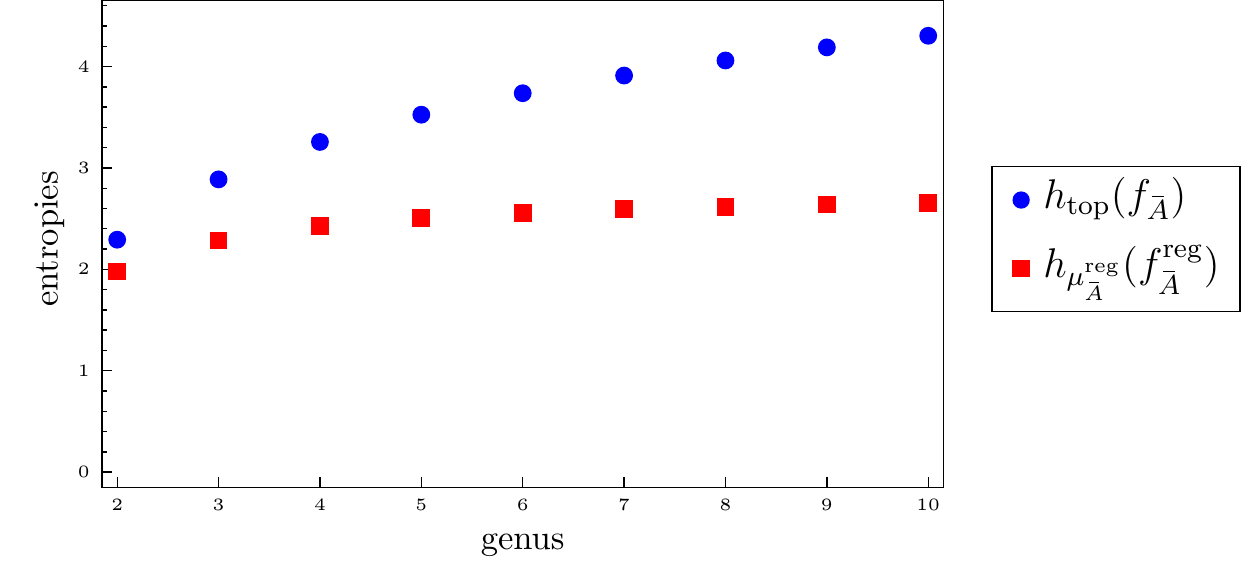}
 \caption{Topological entropy and mea\-sure-theo\-re\-tic entropy for different genera}
 \label{fig entropies}
\end{figure}

\begin{prop}[{\cite[Corollary~7]{AKU-Flexibility}}] \label{thm less}
 If $\A$ is extremal or has short cycles, the mea\-sure-theo\-re\-tic entropy of $f_\A$ with respect to its smooth invariant measure~$\mu_\A$ is strictly less than the topological entropy of $f_\A$. 
\end{prop}

From the expression for $H(g)$ in \Cref{thm max and flexibility}\ref{item max}, the maximum possible value of the measure-theoretic entropy with respect to the smooth invariant measure is 
\[ h_{\mu_\A^\mathrm{reg}}(f_\A^\mathrm{reg}) = \frac{\pi^2(4g-4)}{(8g\!-\!4)\arccosh(1\!+\!2\cos\tfrac\pi{4g-2})}, \]
and we show in~\cite[Section~4]{AKU-Flexibility} that this is strictly larger than the topological entropy 
\[ h_\mathrm{top}(f_\A) = \log\big( 4g - 3 + \sqrt{(4g-3)^2 - 1}\big) \]
for all $g \ge 2$. See \Cref{fig entropies} for a graph of these values.

\subsection{Examples of other Markov partitions} \label{sec sofic examples}

Only countably many multi-para\-me\-ters $\A$ will admit a Markov partition for $f_\A$. Specifically, Markov cases occur when the right and left orbits of $A_k$, that is,
\[\{ \lim_{\varepsilon\to0^+} f_\A^n(A_k+\varepsilon)\}_{n=0}^\infty  \quad\text{and}\quad \{ \lim_{\varepsilon\to0^+} f_\A^n(A_k-\varepsilon)\}_{n=0}^\infty, \]
are eventually periodic for each $\A_k$ (these are called ``upper'' and ``lower'' orbits in~\cite{KU17}). For extremal $\A$, this is immediate. If $\A$ has the cycle property, it is sufficient that the ``cycle ends'' are discontinuity points, as stated for short cycles in~\cite[Proposition 7.2]{AK19}.

\begin{figure}[hbt]
\includegraphics[width=0.9\textwidth]{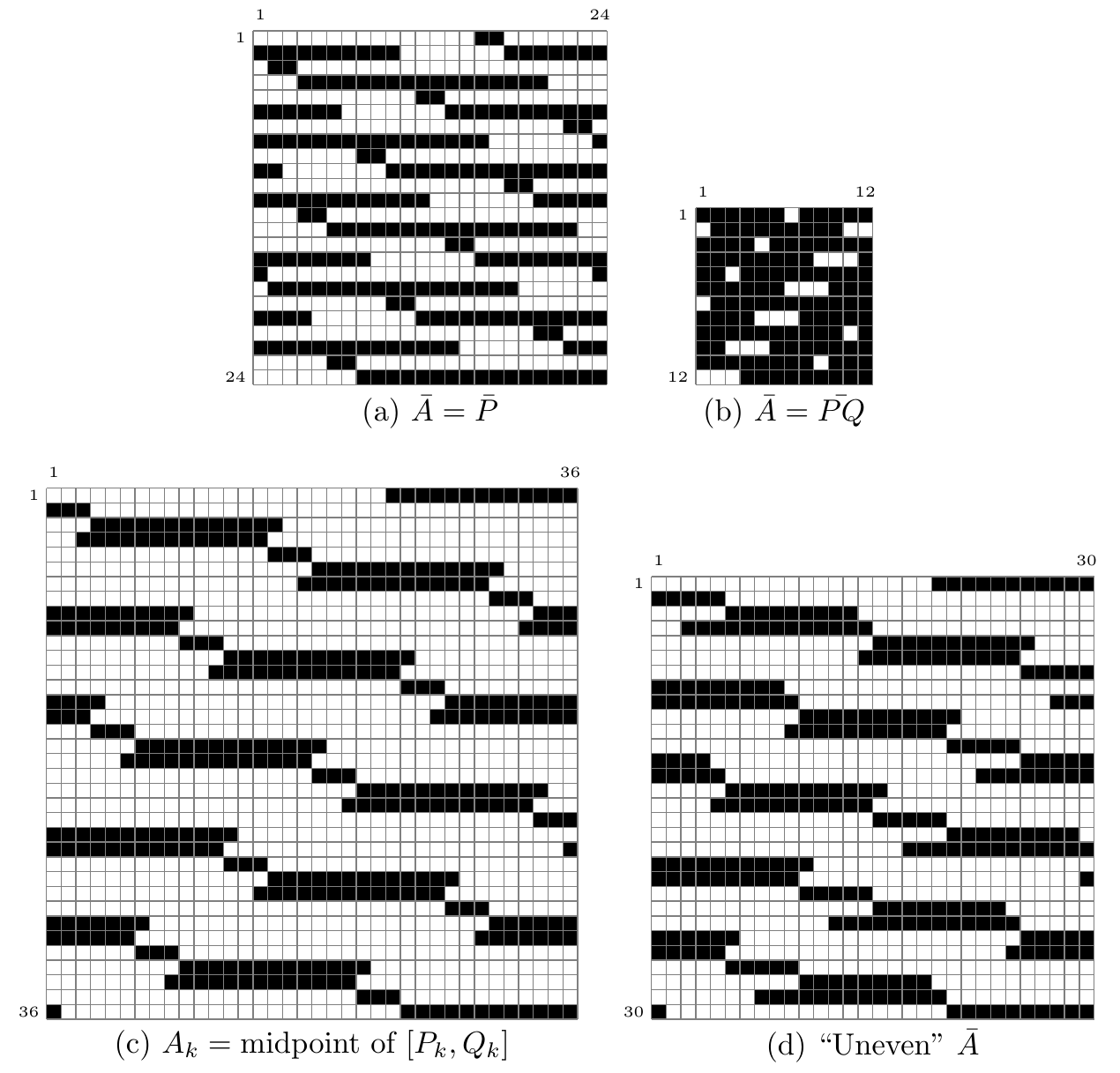}
\caption{Four transition matrices with maximal eigenvalue $5+2\sqrt6$}
\label{fig matrices}
\end{figure}

When $f_\A$ is Markov, the maximal eigenvalue of the associated transition matrix must be $\lambda = 4g-3 + \sqrt{(4g-3)^2-1}$ because we know that $h_\mathrm{top}(f_\A)$ is the $\log$ of this value for all parameters $\A$. Given a Markov partition $\{I_1,...,I_n\}$ with $n$ elements, the entries of the associated $n \times n$ transition matrix are described in \eqref{transition matrix defn}.
In this section, we give four examples of maps and transition matrices with, respectively,
\begin{enumerate}[label=(\alph*)] \addtolength\itemsep{1mm}
	\item \label{item matrix P} $n = 2(8g-4)$
	\item \label{item matrix PQ} $n = 8g-4$
	\item \label{item matrix midpoint} $n = 3(8g-4)$
	\item \label{item matrix uneven} $n = \tfrac52(8g-4)$
\end{enumerate}
partition elements. \Cref{fig matrices} shows these matrices for genus~$2$ (each white cell is a~$0$ in the matrix, and each black cell is a~$1$); each matrix is for a different map~$f_\A$ and uses a different partition~$\{I_1,...,I_n\}$, but all have the same maximal eigenvalue.

\smallskip\textbf{\ref*{item matrix P}} Let $A_k = P_k$ for all $k$. As described in~\cite{AF91,AKU-Flexibility}, the partition $\{I_1,\dots,I_{16g-8}\}$ with
\[ 
	I_{2k-1} := [P_k, Q_k], \quad 
	I_{2k} := [Q_{k}, P_{k+1}], \qquad 
	k=1,\dots,8g-4,
\]
is a Markov partition for $f_\P$ (and, in fact, for every extremal $\A$). For $g=2$, the matrices $M_\P$ and $M_\Q$ for this partition are shown in~\Cref{fig PQ matrices}, and~$M_\P$ is shown again in \Cref{fig matrices}\ref*{item matrix P} as a colored grid.

\smallskip\textbf{\ref*{item matrix PQ}} For $\A = \bar{PQ} = \{P_1,Q_2,P_3,Q_4,...,\}$ we \textit{could} use the partition $\{I_1,\dots,I_{16g-8}\}$ from \ref{item matrix P} and~\cite{AKU-Rigidity}, but the coarser partition $\{I_1$, ..., $I_{8g-4}\}$ with
\[
	I_k = [P_k,Q_{k+1}] \text{ for odd $k$}, \quad 
	I_k = [Q_k,P_{k+1}] \text{ for even $k$}
\]
is also Markov for $f_{\bar{PQ}}$. Using this partition we get a much smaller transition matrix (it is $(8g-4) \times (8g-4)$ rather than $(16g-8) \times (16g-8)$) with exactly the same maximal eigenvalue. For $g=2$, this $12 \times 12$ matrix is shown in \Cref{fig matrices}\ref*{item matrix PQ}. 

\smallskip\textbf{\ref*{item matrix midpoint}} Let $A_k$ be the midpoint of $[P_k,Q_k]$ for each $k$. As described in~\cite[Sections 7-8]{AK19}, this $f_\A$ has a Markov partition $\{I_1,...,I_{24g-12}\}$ with
\begin{align*}
	I_{3k-2} := [A_k, ], \quad 
	I_{3k-1} := [C_k, B_k], \quad
	I_{3k} := [B_k, A_{k+1}], \qquad
	k=1,\dots,8g-4,
\end{align*}
where \[ B_k := T_{\sigma(k-1)}A_{\sigma(k-1)} \quad\text{and}\quad C_k := T_{\sigma(k+1)}A_{\sigma(k+1)+1} \] are the first iterates in the right and left orbits. For $g=2$, the corresponding $36 \times 36$ transition matrix is shown in \Cref{fig matrices}\ref*{item matrix midpoint}.

\smallskip\textbf{\ref*{item matrix uneven}} For odd $k$ let $A_k$ be the midpoint of $[P_k,Q_k]$, and for even $k$ let $A_k$ be the image under $T_{\sigma(k)+1} \circ T_k$ of the midpoint of $[P_{k+1},Q_{k+1}]$. In~\cite[Sec\-tion~8, Example~3]{AK19} we show that this $f_\A$ has a partition with $3(8g-4)$ intervals following the same construction as \ref{item matrix midpoint}. %, so one can construct a $(24g-12) \times (24g-12)$ transition matrix using that partition. 
However, for this particular~$\A$ we have that
$ C_{k+1} = B_{k+1} \text{ if $k$ is odd,} $
which means that the intervals $I_5 = I_{3(2)-1} = [C_2,B_2]$, $I_{11} = I_{3(4)-1} = [C_4,B_4]$, etc., are trivial. These comprise $4g-2$ of the $24g-12$ intervals.

Thus we can make a Markov partition for this $f_\A$ using only $20g-10$ intervals. Explicitly, we use $\{I_1,...,I_{20g-10}\}$ given by
\begin{align*}
	I_{5k-4} := [A_k, C_k], \quad 
	I_{5k-3} := [C_k, B_k], \quad
	I_{5k-2} := [B_k, A_{k+1}], \quad \\*
	I_{5k-1} := [A_{k+1}, B_{k+1}], \quad
	I_{5k} := [B_{k+1}, A_{k+2}], \qquad
	k=1,\dots,4g-2.
\end{align*}
For $g=2$, the corresponding $30 \times 30$ transition matrix is shown in \Cref{fig matrices}\ref*{item matrix uneven}.

\subsection{\texorpdfstring{Parameters for $\boldsymbol{\mathcal T(2)}$}{Teichm\"uller parameters for genus 2}} \label{sec genus 2}

We now provide some geometric results specific to genus $2$. We use the parameterization of~$\Tc(2)$ given by Maskit in~\cite{M99}:

\begin{thm}[{\cite[Theorem~8.1]{M99}}] \label{thm parameterization}
For any $(\alpha,\beta,\gamma,\sigma,\tau,\rho) \in \R^6$, define
\begin{equation} \label{greeks} \begin{split}
	\mu :\hspace{-3pt} 
	&= \arccosh( \coth\beta \cosh\sigma \cosh\tau + \sinh\sigma \sinh\tau ), \\
	\delta :\hspace{-3pt} 
	&= \operatorname{arccoth}\!\bigg( \frac{ \cosh\gamma \cosh\mu - \coth\alpha \sinh\gamma \sinh\mu - \sinh\rho \sinh\sigma }{\cosh\rho \cosh\sigma} \bigg), \\
	\epsilon :\hspace{-3pt} 
	&= \arccosh( \coth\mu \sinh\alpha \sinh\gamma - \cosh\alpha \cosh\gamma ), \\
	\phi :\hspace{-3pt} 
	&= \arccosh\!\bigg( \frac{ \sinh\beta \sinh\delta (\sinh\rho \sinh\tau + \cosh\gamma) }{\cosh\rho \cosh\tau} - \cosh\beta \cosh\delta \bigg).
\end{split} \end{equation}
There exists a real-analytic diffeomorphism between
\[ \big\{ (\alpha,\beta,\gamma,\sigma,\tau,\rho) \in \R^6 : \alpha>0,\beta>0,\gamma>0, \delta > 0, \epsilon > 0 \big\} \]
and the Teichm\"uller space $\Tc(2)$. (Note that $\delta > 0$ and $\epsilon > 0$ are implicitly requirements on $\alpha,\beta,\gamma,\sigma,\tau,\rho$ and that the stated inequalities imply $\phi > 0$.)
\end{thm}

In~\cite{M99}, $\mu$ and $\delta$ are defined geometrically and are then proven to satisfy Equation~\cite[\!(12)]{M99}, that is,
\[ \coth\delta = \frac{\cosh\gamma \cosh\mu - \coth\alpha \sinh\gamma \sinh\mu - \sinh\sigma \sinh\rho}{\cosh\sigma \cosh\rho}, \]
when all parameters satisfy certain necessary inequalities. %\cite[\!(2), (3), (13)]{M99}.
Here, we are defining~$\mu$ and~$\delta$ by~\eqref{greeks} and making $\delta > 0$ a requirement of valid parameters. This condition $\delta > 0$ is equivalent to the inequality~\cite[\!(13)]{M99}. Maskit does not explicitly use $\epsilon$ or $\phi$ at all, but algebraic manipulation shows that~\cite[\!(3)]{M99} is equivalent to $\epsilon > 0$. Additionally, the inclusion of $\epsilon$ and $\phi$ allows for the concise statement of \Cref{thm perimeter}.

\smallskip The correspondence between parameters $(\alpha,\beta,\gamma,\sigma,\tau,\rho)$ and polygons $\Fc$ can be made explicit as follows. Let
\pagebreak 
\providecommand\temp{}
\renewcommand\temp[4]{ X_{#4} \begin{pmatrix} #2 & 0 \\ 0 & #3 \end{pmatrix} X_{#4}^{-1} \qquad\text{with } X_{#4} = \begin{pmatrix}i&1\\1&i\end{pmatrix} \! \begin{pmatrix} #1 \end{pmatrix} }
\begin{align*}
	A &= \temp{ 1 & 1 \\ e^\mu & e^{-\mu} }{ e^\alpha }{ e^{-\alpha} }{A}\!, \\
	B &= \temp{ e^{\sigma} & e^{\sigma} \\ e^{-\tau} & -e^{\tau} }{ e^\beta }{ e^{-\beta} }{B}\!, \\
	C &= X_C \begin{pmatrix} e^\gamma & 0 \\ 0 & e^{-\gamma} \end{pmatrix} X_C^{-1} \qquad\text{with } X_C = \begin{pmatrix}i&1\\1&i\end{pmatrix}\!, \\
	D &= \temp{ e^{\sigma+\gamma} & e^{\sigma+\gamma} \\ -e^{\rho} & e^{-\rho} }{ e^\delta }{ e^{-\delta} }{D}\!,
\end{align*}
$E = A^{-1}C^{-1}$, and $F = D^{-1}B^{-1}$ (the eigenvalues of $E$ and $F$ are, respectively, $e^{\pm\epsilon}$ and $e^{\pm\phi}$, but the matrices that diagonalize $E$ and $F$ do not seem to have simple expressions). Then the sides of the polygon $\Fc$ lie along axes of M\"obius transformations that are products of these matrices; explicitly, if \renewcommand\SS{s}$\SS_k$ is the transformation for which~$P_k$ is the repelling fixed point and~$Q_{k+1}$ is the attracting fixed point, then
\begin{align*}
	\SS_1 &= C^{-1}D^{-1}C, &	
	\SS_2 &= AC, &
	\SS_3 &= AF^{-1}A^{-1}, &
	\SS_4 &= A^{-1}, \\*
	\SS_5 &= F, & 
	\SS_6 &= E^{-1}, &
	\SS_7 &= D, &
	\SS_8 &= DED^{-1}, \\*
	\SS_9 &= B^{-1}D^{-1}, &
	\SS_{10} &= B^{-1}AB, &
	\SS_{11} &= \SS_1^{-1} B, &
	\SS_{12} &= C^{-1} \SS_{10}^{-1}.
\end{align*}
These expressions for $\SS_k$ are from~\cite[Appendix]{AKU-Flexibility} (called ``$S_k$'' there), and the expres\-sions for $A, B, C, D$ are newly-discovered formulas that are easily shown to be equivalent to the original descriptions of the matrices in~\cite[Theorem~5.1]{M99}.

\begin{prop} \label{thm perimeter} 
    Let $(\alpha,\beta,\gamma,\sigma,\tau,\rho)$ be Maskit's parameters for a point in $\Tc(2)$, and let $\Fc$ be the associated fundamental polygon. Then \[ \mathrm{Perimeter}(\Fc) = 4(\alpha\!+\!\delta\!+\!\epsilon\!+\!\phi) \] with $\delta, \epsilon, \phi$ defined by \eqref{greeks}. This immediately implies {$h_{\mu_\A}(f_\A) = {\pi^2}/(\alpha\!+\!\delta\!+\!\epsilon\!+\!\phi)$}. 
\end{prop}

\providecommand\tr{\operatorname{tr}}
\begin{proof}
    For a M\"obius transformation~$M \in \G$, the axis of~$M$ on~$\D$ pro\-jects to a closed geodesic on~$\G\backslash\D$, and the length of this closed geodesic is $2\arccosh|\tr(M)/2|$, where $\tr(M)$ is the trace of~$M$.
    
    Sides $1, 4, 7, 10$ of $\Fc$ are closed geodesics on the surface $S$, while the other sides come in pairs: the closed geodesic along the axis of $\SS_k$ for $k \notin \{1, 4, 7, 10\}$ consists of side $k$ and side $k+6$.
    
    The eigenvalues of $A, B, C, D, E, F$ are $e^{\pm\alpha}, e^{\pm\beta}, e^{\pm\gamma}, e^{\pm\delta}, e^{\pm\epsilon}, e^{\pm\phi}$, respectively.
    Thus the length of side $1$ is 
    \begin{align*}
        2\arccosh|\tr(\SS_1)/2| 
        &= 2\arccosh|\tr(C^{-1}D^{-1}C)/2|
        = 2\arccosh|\tr(D^{-1})/2|
        \\&= 2\arccosh|\tfrac12(e^{\delta} + e^{-\delta})|
        = 2\arccosh|\cosh\delta|
        = 2\delta
    \end{align*}
    and the combined length of sides $2$ and $8$ is
    \begin{align*}
        2\arccosh|\tr(\SS_2)/2| 
        &= 2\arccosh|\tr(AC)/2|
        = 2\arccosh|\tr(CA)/2|
        \\&= 2\arccosh|\tr(E^{-1})/2|
        = 2\arccosh|\tfrac12(e^{\epsilon} + e^{-\epsilon})|
        = 2\epsilon.
    \end{align*}
    Similarly, the combined length of sides $3$ and $9$ is $2\arccosh|\tr(F)/2| = 2\phi$, and the length of side $4$ is $2\arccosh|\tr(A)/2| = 2\alpha$. Sides $6, 7, 10, 11, 12$ have the same length as, respectively, sides $8, 1, 4, 9, 2$ because those respective side pairs are identified. Thus the total perimeter of $\Fc$ is
    \[ 2( 2\delta + 2\epsilon + 2\phi + 2\alpha ) = 4(\alpha+\delta+\epsilon+\phi). \qedhere \]
\end{proof}

\begin{remark} For many surfaces (e.g., all those with $\sigma=\tau=\rho=0$) sides $2$ and $8$ will each have length $\epsilon$, sides $3$ and $9$ will each have length $\phi$, etc., but in general we require only that sides $2$ and $8$ together have total length $2\epsilon$, and so on, with only sides~$1,4,7,10$ having simple expressions for their individual lengths. \end{remark}

\begin{figure}[hbt]
	\includegraphics[width=0.75\textwidth]{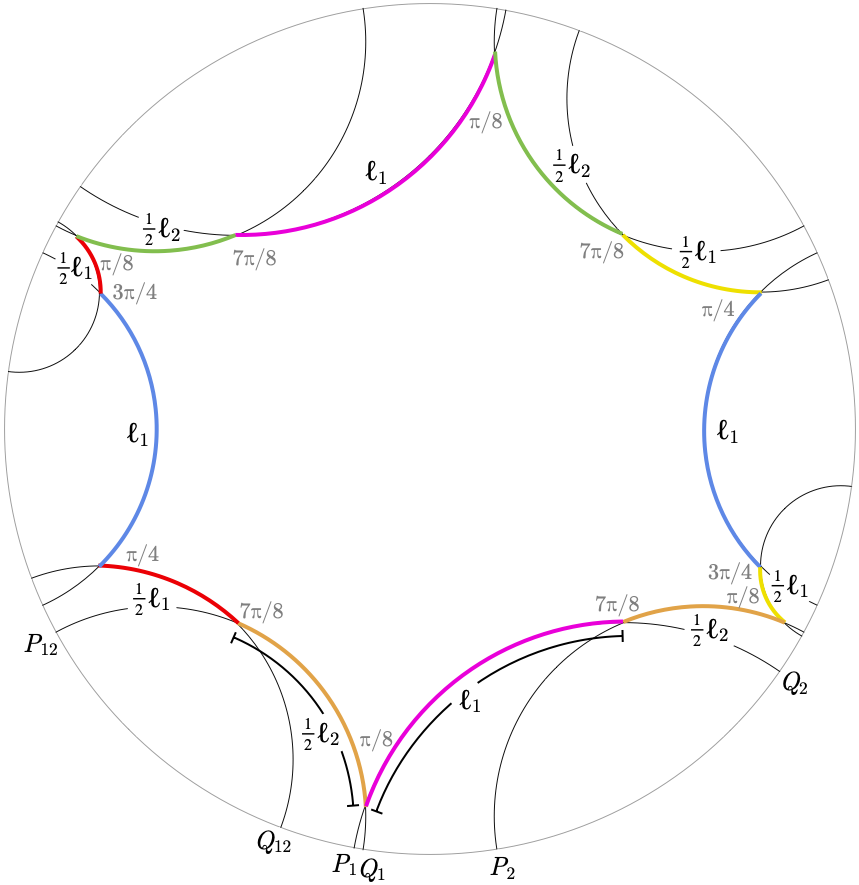}
	\caption{Fundamental $12$-gon for the Bolza surface (sides of the same color are identified)}
	\label{fig bolza}
\end{figure}

\smallskip
The subspace of $\Tc(2)$ in which $\sigma=\tau=0$ contains several notable surfaces (in these cases, the origin of $\D$ is a ``Weierstrass point'' of the surface):
\begin{itemize}
	\item The regular polygon corresponds to $\alpha = \tfrac12 \arccosh(1+\sqrt3)$, $\beta = \gamma = 2\alpha$, and $\rho = 0$. The perimeter of the regular polygon for $g=2$ is $12 \arccosh(1 + \sqrt3) \approx 19.955$, and the associated entropy is $H(2) = \frac{4\pi^2}{12 \arccosh(1 + \sqrt3)} \approx 1.978$.
	\item Maskit's ``base surface''~\cite[Sec.~3]{M99} corresponds to $\alpha = \beta = \gamma = \arccosh(2)$ and $\rho = 0$. The perimeter of $\Fc$ for this surface is $16 \arccosh(2) \approx 21.071$, and so the entropy of the boundary map is ${\approx} 1.874$.
	\item The Bolza surface is the hyperbolic genus $2$ surface that maximizes the systole (the length of the shortest closed geodesic on the surface)~\cite{S93}. Denoting \[ \ell_1 = 2 \arccosh(1 + \sqrt2) \approx 3.057, \qquad \ell_2 = 
 2 \arccosh(3 + 2\sqrt2) \approx 4.897, \] this systole is exactly $\ell_1$.
	This surface corresponds to $\alpha = \beta = \gamma = \tfrac12 \ell_1$ and \mbox{$\rho = \operatorname{arcsinh}(1)$.}\footnote{The authors thank Polina Vytnova for suggesting and assisting with the question of parameters for the Bolza surface.} Usually the Bolza surface is described as a gluing of a regular \emph{octagon}, but it can also be described as a gluing of a $12$-gon, as shown in~\Cref{fig bolza}. Sides 1, 4, 7, 10 have length $\ell_1$; sides 2, 6, 8, 12 have length $\tfrac12 \ell_2$; and sides 3, 5, 9, 11 have length $\tfrac12 \ell_1$. Thus the perimeter of $\Fc$ is $6 \ell_1 + 2 \ell_2 \approx 28.137$, and the entropy of the boundary map is ${\approx} 1.403$. 
\end{itemize}

\section{Open questions} \label{sec open}

The flexibility results of \Cref{sec flexibility} are proved under the assumption that the multi-parameter $\A$ is extremal or has short cycles. In order to say anything about the measure-theoretic entropy $h_{\mu_\A}(f_\A)$, we must have a well-defined measure $\mu_\A$, and this is not guaranteed in all cases. This could be established directly---\Cref{conjecture measure}---or through the connection to geodesic flow and natural extension maps---Conjectures~\ref{conjecture attractor} and~\ref{conjecture conjugacy}. Additionally, we use the conjugacy $\Phi_\A: \Omega_{\rm geo} \to \Omega_\A$ between $F_{\rm geo}$ and $F_\A\big|_{\Omega_\A}$ to prove that $h_{\mu_\A}(f_\A) = \pi \cdot \mathrm{Area}(\Fc) / \nu(\Omega_\A)$ from Step~2 of \Cref{sec flexibility}, and this is currently known only for extremal and short cycle multi-parameters.

\begin{conjecture} \label{conjecture measure}
    For any $\Fc \in \Tc(g)$ and any $\A$ with $A_k \in [P_k,Q_k]$, there exists a smooth $f_\A$-invariant ergodic probability measure $\mu_\A$.
\end{conjecture}

The existence of this measure for maps $f_\A$ admitting a Markov partition can be established directly using Adler's ``Folklore Theorem,'' similarly to the case $\bar A=\P$ discussed in~\cite{AF91, BS79, B79}. For regular fundamental polygons, existence can be established using results about expanding, non-Markov maps in~\cite{Zw}.

\begin{conjecture} \label{conjecture attractor}
    For any $\A$ with $A_k \in [P_k,Q_k]$, there exists a set $\Omega_\A \subset \Sb \times \Sb \setminus \Delta$ with finite rectangular structure that is a domain of bijectivity for $F_\A$ and moreover the global attractor of $F_\A$.
\end{conjecture}

\Cref{conjecture attractor} is part of the ``Reduction Theory'' proposed by Don Zagier and described for Fuchsian groups in \cite[Introduction]{KU17}. Additionally, understanding the structure of $\Omega_\A$ may help in proving the following:

\begin{conjecture} \label{conjecture conjugacy}
    For any $\A$ with $A_k \in [P_k,Q_k]$, the map $F_\A\big|_{\Omega_\A}$ is conjugate to $F_{\rm geo}$ by a map $\Phi_\A:\Omega_{\rm geo} \to \Omega_\A$ that acts piecewise by M\"obius transformations.
\end{conjecture}

If \Cref{conjecture conjugacy} is true then \Cref{thm max and flexibility} in fact holds for all $\A$ with $A_k \in [P_k,Q_k]$.

\end{document}